\documentclass[11pt,twoside,a4paper]{article}

\usepackage{amsmath,amssymb}
\usepackage{graphicx, epsfig, float}
\usepackage{stmaryrd, verbatim}
\usepackage{tikz}
\usepackage{hyperref}
\usepackage{lineno}
\usepackage{authblk}

\newcommand{\R}{ {\mathbb R} }
\newcommand{\bZ}{ {\mathbb Z} }

\newcommand{\ep}{\epsilon}

\newcommand{\Om}{\Omega}
\newcommand{\Omo}{\Omega_{1\ep}}
\newcommand{\Omt}{\Omega_{2\ep}}

\newcommand{\strain}[1]{ e ( {\pmb #1} ) }

\newcommand{\dv}[1]{{\rm div} \ #1}
\newcommand{\curl}[1]{{\rm curl} \ #1 }

\newcommand{\jump}[1]{\left \llbracket #1 \right \rrbracket} 
\newcommand{\vc}[1]{{\pmb #1}}

\newcommand{\grad}[1]{\nabla{#1}}

\newtheorem{remark}{Remark}

\newcommand{\subf}[2]{%
  {\small\begin{tabular}[t]{@{}c@{}}
  #1\\#2
  \end{tabular}}%
}

\begin{document}

\title{Multiscale modeling of magnetorheological suspensions}
\author[1]{Grigor Nika\footnote{Corresponding author: \url{grigor.nika@wias-berlin.de}}}
\author[2]{Bogdan Vernescu}
\affil[1]{Weierstrass Institute for Applied Analysis and Stochastics\\ Mohrenstra\ss{}e 39, 10117 Berlin, Germany}
\affil[2]{Department of Mathematical Sciences\\ Worcester Polytechnic Institute\\ 100 Institute Rd.\\ Worcester, MA 01601}

\maketitle

\begin{abstract}
\noindent
We develop a multiscale approach to  describe the behavior of a suspension of solid magnetizable particles in a viscous non-conducting fluid in the presence of an externally applied magnetic field. By upscaling the quasi-static Maxwell equations coupled with the Stokes' equations we are able to capture the magnetorheological effect. The model we obtain generalizes the one introduced by Neuringer \& Rosensweig~\cite{NR64}, \cite{Ro14} for quasistatic phenomena. We derive the macroscopic constitutive properties explicitly in terms of the solutions of local problems. The effective coefficients have a nonlinear dependence on the volume fraction when {\it chain structures} are present. The velocity profiles computed for some simple flows, exhibit an {\it apparent yield stress} and the flow profile resembles a {\it Bingham fluid} flow.
\end{abstract}


\section{Introduction}
Magnetorhelogical fluids are a suspension of non--colloidal, ferromagnetic particles in a non--magnetizable carrier fluid. The particles are often of micron size ranging anywhere from $0.05-10\, \mu m$ with particle volume fraction from $10-40 \, \%$. They were discovered by J. Rabinow in $1948$ \cite{Ra48}. Around the same time W. Winslow~\cite{Win49} discovered electrorheological fluids, a closely related counterpart~\cite{Bossis}, \cite{HL11}, \cite{PV02}, \cite{Ver02}.  

Magnetorheological fluids respond to an external magnetic field by a rapid, reversible change in their properties. They can transform from a liquid to a semi solid state in a matter of milliseconds. Upon the application of a magnetic field, the dipole interaction of adjacent particles aligns the particles in the direction of the magnetic field lines. Namely particles attract one another along the magnetic field lines and repel one another in the direction perpendicular to them. This leads to the formation of {\it aggregate structures} or {\it chain structures}. Once these {\it chain structures} are formed, the magnetorheological fluid exhibits a higher viscosity and yield stress that can now be triggered and dynamically controlled by an applied external magnetic field~\cite{Bossis}, \cite{KD17}, \cite{Ha92}, \cite{Lo06}.

The presence of these {\it chain structures} leads to a non-Newtonian behavior of the fluid. In many works, the Bingham constitutive law is used as an approximation to model the response of the magnetorheological and electrorheological fluids, particularly in shear experiments~\cite{GS15}, \cite{HL11}, \cite{PV02}, \cite{CP17}. Although the Bingham model has proven itself useful in characterizing the behavior of magnetorheological fluids, it is not always sufficient. Recent experimental data show that true magnetorheological fluids exhibit departures from the Bingham model~\cite{HL11}, \cite{YLL05}, \cite{ZFF17}.

Another member of the magnetic suspensions family are ferrofluids. Ferrofluids are stable colloidal suspensions of nanoparticles in a non-magnetizable carrier fluid. The initiation into the hydrodynamics of ferrofluids began with Neuringer and Rosensweig in $1964$ \cite{NR64} and by a series of works by Rosensweig and co-workers summarized in \cite{Ro14}. The model introduced in \cite{NR64} assumes that the magnetization is collinear with the magnetic field and has been very useful in describing quasi-stationary phenomena. This work was extended by Shliomis \cite{Sh72} by avoiding the collinearity assumption of the magnetization and the magnetic field and by considering the rotation of the nanoparticles with respect to the fluid they are suspended in.

Most of the models characterizing magnetorheological suspensions are derived phenomenologically. The first attempt to use homogenization mechanics to describe the behavior of magnetorheological / electrorheological fluids was carried out in \cite{Le85}, \cite{LH88} and \cite{PV02}. In the works \cite{Le85}, \cite{LH88} the influence of the external magnetic field is introduced as a volumic density force acting on each particle and as a surface density force acting on the boundary of each particle. The authors in \cite{PV02} extend the work in \cite{LH88}, for electrorheological fluids, by presenting a more complete model that one way couples the conservation of mass and momentum equations with Maxwell's equations through the {\it Maxwell stress tensor}. As an application they consider a uniform shearing of the electrorheological fluid submitted to a uniform electric field boundary conditions in a two dimensional slab and they recover a stress tensor, at the macroscopic scale, that has exactly the form of the Bingham constitutive equation.

The authors in~\cite{PV02}, \cite{Ro14}, \cite{Ru02}, \cite{Ver02}, \cite{HMF16} use models that decouple the conservation of mass and momentum equations from the Maxwell equations. Thus, in principle, one can solve the Maxwell equations and use the resulting magnetic or electric field as a force in the conservation of mass and momentum equations. 

The present work focuses on a suspension of rigid magnetizable particles in a Newtonian viscous fluid with an applied external magnetic field. We assume the fluid to be electrically non-conducting. We use the homogenization method to upscale the quasi-static Maxwell equations coupled with the Stokes equations through Ohm's law to capture the magnetorheological effect. In doing so we extend the model of \cite{PV02}, \cite{Ver02}. Thus, the Maxwell equations, and the balance of mass and momentum equations must be simultaneously solved. Additionally, the model is able to capture the added effect particle {\it chain structures} have on the effective coefficients. We demonstrate this added effect by carrying out explicit computations of the effective coefficients using the finite element method.   

The paper is organized in the following way. In Section 2. we introduce the problem in the periodic homogenization framework. The particles are periodically distributed and the size of the period is of the same order as the characteristic length of the particles. We assume the fluid velocity is continuous across the particle interface and that the particles are in equilibrium in the presence of the magnetic field.

In Section 3. we use two-scale expansions to obtain a one way coupled set of local problems at order $\mathcal{O}(\ep^{-1})$. One problem characterizes the effective viscosity of the magnetorheological fluid while the other local problem represents the magnetic field contribution.

Section 4. and Section 5. are devoted to the study of the local problems that arise from the contribution of the bulk magnetic field as well as the bulk velocity and we provide new constitutive laws for Maxwell's equations. 

In Section 6. we provide the governing effective equations of the magnetorheological fluid which include, in addition to the viscous stresses, a ``{\it Maxwell type}'' stress of second order in the magnetic field. Furthermore, we provide formulas for the effective viscosity and three different effective magnetic permeabilities for the ``{\it Maxwell type}" stress that generalize those in \cite{LH88}.

Section 7. is devoted to numerical results for a suspension of circular iron particles of different volume fractions using the finite element method. Moreover, we explore the effect that {\it chain structures} have in the effective coefficients and compare the results with the absence of {\it chain structures}. Moreover, we compute the velocity profiles for Poiseuille and Couette flows of the magnetorheological fluid and we plot the shear stress versus the shear rate curve for different values of the applied magnetic field to obtain a yield stress comparable to the one observed experimentally (e.g. \cite{YLL05}). Finally Section 8. contains conclusions and perspectives on the study.
 
\subsection*{Notation}
Throughout the paper we are going to be using the following notation: $I$ indicates the $n \times n$ identity matrix, bold symbols indicate vectors in two or three dimensions, regular symbols indicate tensors, $e(\vc{u})$ indicates the strain rate tensor defined by $\displaystyle e(\vc{u}) = \frac{1}{2} \left ( \nabla \vc{u} + \nabla \vc{u}^\top \right)$, where often times we will use subscript to indicate the variable of differentiation. The inner product between matrices is denoted by  $A$:$B$ = $tr(A^\top B) = \sum_{ij} A_{ij} \, B_{ji}$ and throughout the paper we employ the Einstein summation notation for repeated indices.

\section{Problem statement}
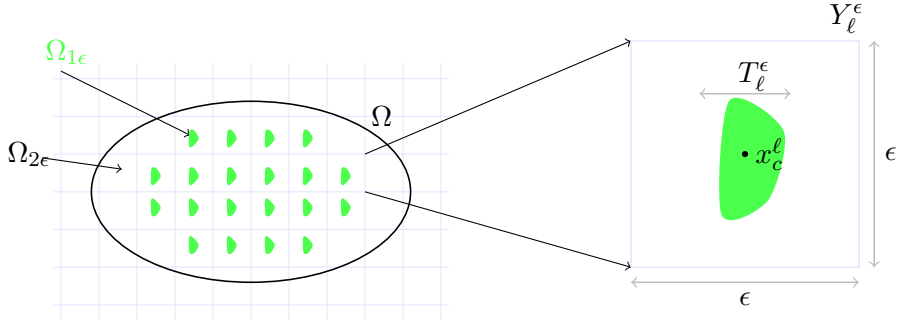
\begin{figure}[!htb]
\begin{center}
\begin{tikzpicture}[scale=1.0]
% Draw thin grid lines with color 40% green!70 + 60% white
\draw [step=0.5,thin,blue!10] (-2.6,-1.7) grid (2.6,1.7);

% Draw x and y axis lines
%\draw [->] (-3.5,0) -- (3.5,0) node [below] {$x$};
%\draw [->] (0,-2.5) -- (0,2.5) node [left] {$y$};

%Drawing of the domain. It is not necessarilly and ellipse but this is the best shape I could find.
\draw [semithick,black] (0,0) ellipse (2.1 and 1.2);

\draw [semithick,black] (2.0,1.0) node [left] {$\Om$};

\draw [semithick,black] (-2.0,1.85) node [left] {$\color{green!70} \Omo$};
\draw [semithick,black] (-2.5,0.5) node [left] {$\Omt$};

\draw [->] (-2.5,1.6) -- (-0.8,0.75); % node [below] {$x$};
\draw [->] (-2.75,0.45) -- (-1.7,0.3); % node [below] {$x$};

\draw[semithick,green!70,fill=green!70] plot [smooth cycle] coordinates {(0.25,0.3) (0.3,0.2) (0.2,0.1) (0.2,0.3)}; 
\draw[semithick,green!70,fill=green!70] plot [smooth cycle] coordinates {(0.75,0.3) (0.8,0.2) (0.7,0.1) (0.7,0.3)};
\draw[semithick,green!70,fill=green!70] plot [smooth cycle] coordinates {(1.25,0.3) (1.3,0.2) (1.2,0.1) (1.2,0.3)};
%\draw[semithick,green!70,fill=green!70] plot [smooth cycle] coordinates {(1.7,0.2) (1.75,0.3) (1.7,0.3) (1.7,0.2)};

\draw[semithick,green!70,fill=green!70] plot [smooth cycle] coordinates {(0.25,0.8) (0.3,0.7) (0.2,0.6) (0.2,0.8)};
\draw[semithick,green!70,fill=green!70] plot [smooth cycle] coordinates {(0.75,0.8) (0.8,0.7) (0.7,0.6) (0.7,0.8)};
%\draw[semithick,green!70,fill=green!70] plot [smooth cycle] coordinates {(1.2,0.73) (1.25,0.81) (1.2,0.77) (1.2,0.65)};

\draw[semithick,green!70,fill=green!70] plot [smooth cycle] coordinates {(-0.25,0.3) (-0.2,0.2) (-0.3,0.1) (-0.3,0.3)}; 
\draw[semithick,green!70,fill=green!70] plot [smooth cycle] coordinates {(-0.75,0.3) (-0.7,0.2) (-0.8,0.1) (-0.8,0.3)};
\draw[semithick,green!70,fill=green!70] plot [smooth cycle] coordinates {(-1.25,0.3) (-1.2,0.2) (-1.3,0.1) (-1.3,0.3)};
%\draw[semithick,green!70,fill=green!70] plot [smooth cycle] coordinates {(1.7,0.2) (1.75,0.3) (1.7,0.3) (1.7,0.2)};

\draw[semithick,green!70,fill=green!70] plot [smooth cycle] coordinates {(-0.25,0.8) (-0.2,0.7) (-0.3,0.6) (-0.3,0.8)};
\draw[semithick,green!70,fill=green!70] plot [smooth cycle] coordinates {(-0.75,0.8) (-0.7,0.7) (-0.8,0.6) (-0.8,0.8)};
%\draw[semithick,green!70,fill=green!70] plot [smooth cycle] coordinates {(-1.2,0.73) (-1.25,0.81) (-1.2,0.77) (-1.2,0.65)};

\draw[semithick,green!70,fill=green!70] plot [smooth cycle] coordinates {(-0.25,-0.3) (-0.2,-0.2) (-0.3,-0.1) (-0.3,-0.3)}; 
\draw[semithick,green!70,fill=green!70] plot [smooth cycle] coordinates {(-0.75,-0.3) (-0.7,-0.2) (-0.8,-0.1) (-0.8,-0.3)};
\draw[semithick,green!70,fill=green!70] plot [smooth cycle] coordinates {(-1.25,-0.3) (-1.2,-0.2) (-1.3,-0.1) (-1.3,-0.3)};
%\draw[semithick,green!70,fill=green!70] plot [smooth cycle] coordinates {(1.7,0.2) (1.75,0.3) (1.7,0.3) (1.7,0.2)};

\draw[semithick,green!70,fill=green!70] plot [smooth cycle] coordinates {(-0.25,-0.8) (-0.2,-0.7) (-0.3,-0.6) (-0.3,-0.8)};
\draw[semithick,green!70,fill=green!70] plot [smooth cycle] coordinates {(-0.75,-0.8) (-0.7,-0.7) (-0.8,-0.6) (-0.8,-0.8)};
%\draw[semithick,green!70,fill=green!70] plot [smooth cycle] coordinates {(-1.2,0.73) (-1.25,0.81) (-1.2,0.77) (-1.2,0.65)};

\draw[semithick,green!70,fill=green!70] plot [smooth cycle] coordinates {(0.25,-0.3) (0.3,-0.2) (0.2,-0.1) (0.2,-0.3)}; 
\draw[semithick,green!70,fill=green!70] plot [smooth cycle] coordinates {(0.75,-0.3) (0.8,-0.2) (0.7,-0.1) (0.7,-0.3)};
\draw[semithick,green!70,fill=green!70] plot [smooth cycle] coordinates {(1.25,-0.3) (1.3,-0.2) (1.2,-0.1) (1.2,-0.3)};
%\draw[semithick,green!70,fill=green!70] plot [smooth cycle] coordinates {(1.7,0.2) (1.75,0.3) (1.7,0.3) (1.7,0.2)};

\draw[semithick,green!70,fill=green!70] plot [smooth cycle] coordinates {(0.25,-0.8) (0.3,-0.7) (0.2,-0.6) (0.2,-0.8)};
\draw[semithick,green!70,fill=green!70] plot [smooth cycle] coordinates {(0.75,-0.8) (0.8,-0.7) (0.7,-0.6) (0.7,-0.8)};
%\draw[semithick,green!70,fill=green!70] plot [smooth cycle] coordinates {(1.2,0.73) (1.25,0.81) (1.2,0.77) (1.2,0.65)};

%Arrows that point to expanded box
\draw [->] (1.5,0) -- (5,-1); % node [below] {$x$};
\draw [->] (1.5,0.5) -- (5,2); % node [below] {$x$};

%expanded box
\draw [semithick,blue!10] (5,-1) -- (8,-1) -- (8,2) -- (5,2) -- (5,-1);

\draw[semithick,green!70,fill=green!70] plot [smooth cycle] coordinates {(6.2,-0.3) (6.8,-0.1) (7,0.8) (6.3,1.2)};
 
%Do not remove line below. It needs to come after the "hole" otherwise it will not show.
\draw [semithick,black,fill=black] (6.5,0.5) circle (0.03) node [right] {$x^\ell_c$};
\draw [<->,semithick,lightgray] (5.9,1.3) -- (7.1,1.3);
\draw [semithick,black] (6.62,1.20) node [above] {$T^\ep_\ell$};

\draw [semithick,black] (8.2,2.3) node [left] {$Y^\ep_\ell$};

\draw [<->,semithick,lightgray] (8.2,-1) -- (8.2,2);
\draw [semithick,black] (8.2,0.5) node [right] {$\ep$};

\draw [<->,semithick,lightgray] (5,-1.2) -- (8,-1.2);
\draw [semithick,black] (6.5,-1.2) node [below] {$\ep$};

% Draw a triangle with vertices (0,0), (1,0), (1,0.7)
%\draw [semithick,blue] (0,0) -- (1,0) -- (1,0.7) -- cycle;
\end{tikzpicture}
\end{center}
\caption{Schematic of the periodic suspension of rigid magnetizable particles in a non-magnetizable fluid. $\Omo$ represents the domain occupied by the rigid particles and $\Omt$ represents the domain occupied by the surrounding fluid. $\vc{x}_c^\ell$ represents the center of mass of the particle inside the cell of size $\ep$, $Y_\ell^\ep$.}
\label{fig:hom_schem}
\end{figure}

For the homogenization setting of the suspension problem we define $\Om \subset \R^n, \, n=2,3$, to be a bounded open set with sufficiently smooth boundary $\partial \Om$, $Y=\left(-\dfrac{1}{2},\dfrac{1}{2}\right)^n$ is the unit cube in $\R^n$, and $\bZ^n$ is the set of all $n$--dimensional vectors with integer components. 
For every positive $\ep$, let $N^\ep$ be the set of all points $\ell\in\bZ^n$ such that $\ep(\ell+Y)$ is 
strictly included in $\Om$ and denote by $|N^\ep|$ their total number. Let $T$ be the closure of an open 
connected set with sufficiently smooth boundary, compactly included in $Y$. For every $\ep > 0$ and $\ell\in N^\ep$ 
we consider the set $T^\ep_\ell\subset\subset \ep(\ell+Y)$, where $T^\ep_\ell=\ep (\ell+ T)$. The set $T^\ep_\ell$ represents one of the rigid particles suspended in the fluid, and $S^\ep_\ell = \partial T^\ep_\ell$ denotes its surface (see Figure~\ref{fig:hom_schem}). Based on the above setting we define the following subsets of $\Om$:

\[ \Omo=\displaystyle\bigcup_{\ell\in N^\ep} T^\ep_\ell\ ,\quad \Omt= \Om \backslash \overline{\Om}_{1\ep}. \]

In what follows $T^\ep_\ell$ will represent the magnetizable rigid particles, $\Omo$ is the domain occupied by the rigid particles and $\Omt$ the domain occupied by the surrounding fluid. By $\vc{n}$ we indicate the unit normal on the particle surface pointing outwards and by $\jump{ \cdot}$ we indicate the jump discontinuity between the fluid and the rigid part. 

We consider the Navier-Stokes equations coupled with the quasistastic Maxwell equations,

\begin{subequations}\label{one}
\begin{align}
\rho \, \frac{\partial \, \vc{v^\ep}}{\partial \, t} + \rho \, (\vc{v^\ep} \cdot \grad{})\vc{v^\ep} &- \dv{ \sigma^\ep } = \rho \, \vc{f}, \text{ where } \sigma^\ep = 2 \, \nu \, \strain{v^\ep} - p^\ep I & \text{in } &\Omt, \label{st0} \\ 
\dv{\vc{v}^\ep} &= 0, \quad \dv{ \vc{B}^\ep } = 0, \quad \curl{\vc{H}^\ep} = \vc{0} & \text{in } &\Omt, \label{CH1} \\
e(\vc{v}^\ep) &= 0, \quad \dv{ \vc{B}^\ep } = 0, \quad \curl{\vc{H}^\ep} = \eta \, \vc{v^\ep} \times \, \vc{B^\ep} & \text{in } &\Omo, \label{CH2} 
\end{align} 
\end{subequations}
where $\vc{B}^\ep = \mu^\ep \, \vc{H^\ep}$, with boundary conditions on the surface of each particle $T^\ep_\ell$,
\begin{equation}\label{bc1}
\jump{\mathbf{\vc{v}^\ep} } = \vc{0}, \quad \jump{ \vc{B^\ep} \cdot \vc{n} } = 0, \quad \jump{ \vc{n} \times \vc{H^\ep} } = \vc{0} \quad \text{ on $S_\ell^\ep,$ }
\end{equation}
and outer boundary conditions
\begin{equation}\label{bc2}
\vc{v}^{\ep} = \vc{0}, \quad \vc{H}^{\ep} = \vc{b} \quad \text{ on } \partial \Om,
\end{equation}
where $\rho$ is the density of the fluid, $\nu$ is the viscosity, $\vc{v}^\ep$ represents the velocity field, $p^\ep$ the pressure, $\strain{v^\ep}$ the strain rate, $\vc{f}$ the body forces, $\vc{n}$ the exterior normal to the particles, $\vc{H^\ep}$  the magnetic field, $\mu^\ep$ is the magnetic permeability of the material, $\mu^\ep(\vc{x})=\mu_1$ if $ \vc{x} \in \Om_{1\ep}$ and $\mu^\ep(\vc{x})=\mu_2$ if $ \vc{x} \in \Om_{2\ep}$, $\eta$ the electric conductivity of the rigid particles, and $\vc{b}$ is an applied constant magnetic field on the exterior boundary of the domain $\Om$. 

In order to obtain the balance of forces and torques for the particles, let us observe when the magnetorheological fluid is submitted to a magnetic field, the rigid particles are subjected to a force that makes them behave like a dipole aligned in the direction of the magnetic field. This force can be written in the form,
\[
\vc{F^\ep} = -\frac{1}{2} \, | \vc{H^\ep} |^2 \, \grad{\mu^\ep},   
\]
where $| \cdot |$ represents the standard Euclidean norm. The force can be written in terms of the Maxwell stress $\tau_{ij}^\ep = \mu^\ep \, H^\ep_i\, H^\ep_j - \frac{1}{2} \, \mu^\ep \, H^\ep_k\, H^\ep_k\, \delta_{ij}$  as
$\vc{F^\ep} = \dv{\tau^\ep} + \vc{B^\ep} \times \curl{\vc{H^\ep}}.$ Since the magnetic permeability is considered constant in each phase, it follows that the force is zero in each phase. Therefore, we deduce that
\begin{equation}\label{mstress_div}
\dv{\tau^\ep} =
\begin{cases} 
   0 & \text{if } \vc{x} \in \Om_{2\ep} \\
   -\vc{B^\ep} \times \curl{\vc{H^\ep}} & \text{if } \vc{x} \in \Om_{1\ep}.
  \end{cases}
\end{equation}

Lastly, we remark that unlike the viscous stress $\sigma^\ep$, the Maxwell stress is present in the entire domain $\Om$. Hence, we can write the balance of forces and torques in each particle as, 
\begin{gather}\label{bal_forces1}
\begin{aligned}
\int_{T^\ep_\ell} \rho \frac{d\vc{u}^\ep}{dt} \, d\vc{x} &= \int_{S^\ep_\ell} (\sigma^\ep \vc{n} + \jump{\tau^\ep \vc{n}}) \, ds + \int_{T^\ep_\ell} \vc{B^\ep} \times \curl{\vc{H}^\ep} \, d\vc{x} + \int_{T^\ep_\ell} \rho \, \vc{f}\, d\vc{x}, \\
\int_{T^\ep_\ell} \rho (\vc{x}-\vc{x}^{\ell}_c) \times & \frac{d\vc{u}^\ep}{dt} \, d\vc{x} = \int_{S^\ep_\ell} (\sigma^\ep \vc{n} + \jump{\tau^\ep \vc{n}}) \times (\vc{x}-\vc{x^\ell_c}) \, ds \\
   &+ \int_{T^\ep_\ell} (\vc{B^\ep} \times \curl{\vc{H}^\ep}) \times (\vc{x}-\vc{x^\ell_c}) \, d\vc{x} + \int_{T^\ep_\ell} \rho \, \vc{f} \times (\vc{x}-\vc{x^\ell_c})\, d\vc{x},
\end{aligned}
\end{gather}
where $\vc{x_c^\ell}$ is the center of mass of the rigid particle $T_\ell^\ep.$

Equations~\eqref{one}, \eqref{bal_forces1} together with boundary conditions~\eqref{bc1}, \eqref{bc2} describe the behavior of the magnetorheological suspension. 

\subsection{Dimensional Analysis}
Before we proceed further we non-dimensionalize the problem. Denote by $t^*=t/\frac{L}{V}$, $x^* = x / L$, $\vc{v}^* = \vc{v} / V$, $p^* = p / \nu \frac{V}{L}$, $\vc{H}^* = \vc{H} / H$, $\vc{f}^* = \vc{f} / \frac{V^2}{L}$, and $\mu^{\ep*} = \mu^\ep / \mu_2$. Here $L$ is a characteristic length, $V$ is a characteristic velocity, $p$ is a characteristic pressure, $\vc{f}$ is a characteristic force and $H$ is a characteristic unit of the magnetic field. Substituting the above expressions into $\eqref{one}$ as well as in the balance of forces and torques, and using the fact that the flow is assumed to be at low Reynolds numbers, we obtain 
\begin{align*}
Re \left( \frac{\partial \vc{v^\ep}^*}{\partial t} + (\vc{v^\ep}^* \cdot \grad{})\vc{v^\ep}^* \right ) & - \dv^*{ \sigma^{\ep*} } = Re \vc{f}^*, \text{ where } \sigma^{\ep*} = 2 \strain{v^{\ep*}} - p^{\ep*} I & \text{in } &\Omt, \\ 
\dv^*{\vc{v^\ep}^*} &= 0, \quad \dv^*{ \vc{B^\ep}^* } = 0, \quad \curl^*{\vc{H^\ep}^*} = \vc{0} & \text{in } &\Omt,\\
e^*(\vc{v^\ep}^*) &= 0, \quad \dv^*{ \vc{B^\ep}^* } = 0, \quad \curl{\vc{H^\ep}^*} = R_m \vc{v^\ep}^* \times \, \vc{B^\ep}^* & \text{in } &\Omo,
\end{align*} 
where $\vc{B}^{\ep*} = \mu^{\ep*} \, \vc{H^{\ep*}}$ and with boundary conditions on the surface of each particle $T^\ep_\ell$,

\begin{gather*}
\begin{aligned}
\jump{\mathbf{\vc{v}^{\ep*}} } &= \vc{0}, \quad \jump{ \vc{B^{\ep*}} \cdot \vc{n} } = 0, \quad \jump{ \vc{n} \times \vc{H^{\ep*}} } = \vc{0} \quad & \text{ on $S_\ell^\ep,$ }\\
\vc{v}^{\ep*} &= \vc{0}, \quad \vc{H}^{\ep*} = \vc{b^*} &\quad \text{ on } \partial \Om.
\end{aligned}
\end{gather*}
together with the balance of forces and torques,

\begin{gather*}
\begin{aligned}
&\text{Re} \int_{T^\ep_\ell} \frac{d\vc{u}^{\ep*}}{dt^*} d\vc{x}^* = \int_{S^\ep_\ell} \sigma^{\ep*} \vc{n} ds^* + \alpha \int_{S^\ep_\ell} \jump{\tau^{\ep*} \vc{n}} ds^* + \alpha \int_{T^\ep_\ell} \vc{B^\ep}^* \times \curl{\vc{H^\ep}^*} d\vc{x}^*\\
&+\text{Re} \int_{T^\ep_\ell} \vc{f}^* d\vc{x}^*, \\
&\text{Re}\int_{T^\ep_\ell} (\vc{x}^*-\vc{x}^{\ell*}_c) \times \frac{d\vc{u}^{\ep*}}{dt^*} d\vc{x}^* = \int_{S^\ep_\ell} \sigma^{\ep*} \vc{n}\times (\vc{x}^* - \vc{x^\ell_c}^*) ds^* + \alpha \int_{S^\ep_\ell} \jump{\tau^{\ep*} \vc{n}} \times (\vc{x}^*-\vc{x^\ell_c}^*) ds^* \\
&+\alpha \int_{T^\ep_\ell} (\vc{B^\ep}^* \times \curl{\vc{H^\ep}^*}) \times (\vc{x}^*-\vc{x^\ell_c}^*) d\vc{x}^* + \text{Re} \int_{T^\ep_\ell} \vc{f}^* \times (\vc{x}^*-\vc{x^\ell_c}^*) d\vc{x}^*,
\end{aligned}
\end{gather*}
where $\text{Re} = \displaystyle \frac{\rho \, V \, L}{\nu}$ is the Reynolds number, $\alpha = \displaystyle \frac{\mu_2 \, H^2\,L}{\nu\,V}$ is the Alfven number, and $R_m =\displaystyle \eta \, \mu_2 \, L \, V$ is the magnetic Reynolds number. 

In what follows we drop the star for simplicity. Moreover, for low Reynolds numbers the preceding equations become,
\begin{subequations}\label{five}
\begin{align}
-&\dv{ \sigma^{\ep} } = \vc{0}, \text{ where } \sigma^{\ep} = 2 \, \strain{v^{\ep}} - p^{\ep} I &\text{in } &\Omt, \label{five_st0} \\ 
&\dv{\vc{v^\ep}} = 0, \quad \dv{ \vc{H^\ep} } = 0, \quad \curl{\vc{H^\ep}} = \vc{0} & \text{in } &\Omt, \label{five_CH1} \\
&e(\vc{v^\ep}) = 0, \quad \dv{ \vc{H^\ep} } = 0, \quad \curl{\vc{H^\ep}} = R_m \, \vc{v^\ep} \times \, \vc{B^\ep} & \text{in } &\Omo, \label{six_CH2} 
\end{align} 
\end{subequations}
with boundary conditions
\begin{gather}\label{bc_five}
\begin{aligned}
\jump{\mathbf{\vc{v}^\ep} } &= \vc{0}, \quad \jump{ \vc{B^\ep} \cdot \vc{n} } = 0, \quad \jump{ \vc{n} \times \vc{H^\ep} } = \vc{0} \quad & \text{ on $S_\ell^\ep,$ }\\
\vc{v}^{\ep} &= \vc{0}, \quad \vc{H}^{\ep} = \vc{b} &\quad \text{ on } \partial \Om,
\end{aligned}
\end{gather}
together with the balance of forces and torques,
\begin{gather}\label{bal_forces5}
\begin{aligned}
0 &= \int_{S^\ep_\ell} \sigma^{\ep} \vc{n} \, ds + \alpha \, \int_{S^\ep_\ell} \jump{\tau^{\ep} \vc{n}} \, ds + \alpha \, \int_{T^\ep_\ell} \vc{B^\ep} \times \curl{\vc{H^\ep}} \, d\vc{x}, \\
0 &= \int_{S^\ep_\ell} \sigma^{\ep} \vc{n}\times (\vc{x} - \vc{x^\ell_c}) \, ds + \alpha \, \int_{S^\ep_\ell} \jump{\tau^{\ep} \vc{n}} \times (\vc{x} - \vc{x^\ell_c}) \, ds \\
& + \alpha \, \int_{T^\ep_\ell} (\vc{B^\ep} \times \curl{\vc{H^\ep}}) \times (\vc{x} - \vc{x^\ell_c}) \, d\vc{x}.
\end{aligned}
\end{gather}
%In the next section we will use a two scale expansion on the velocity, pressure and the magnetic field. 

\section{Two scale expansions}
We assume the particles are periodically distributed in $\Om$ and thus consider the two scale expansion on $\vc{v^\ep}$, $\vc{H^\ep}$ and $p^\ep$~\cite{All}, \cite{BP89}, \cite{CioDo}, \cite{MV10}, \cite{SP80}, \cite{SPZ87},
\[
	\vc{v^\ep}(\vc{x}) = \sum_{i = 0}^{+\infty} \ep^i \, \vc{v}^i (\vc{x},\vc{y}), \quad \vc{H^\ep}(\vc{x}) = \sum_{i = 0}^{+\infty} \ep^i \, \vc{H}^i (\vc{x},\vc{y}), \quad p^\ep (\vc{x}) = \sum_{i = 0}^{+\infty} \ep^i \, p^i (\vc{x},\vc{y}) \text{ with } \vc{y} = \frac{\vc{x}}{\ep}.
\]
where $\vc{x} \in \Om$ and $\vc{y} \in \R^n.$ One can show that $\vc{v^0}$ is independent of $\vc{y}$ and can thus obtain the following problem at order $\ep^{-1}$,
\begin{subequations}\label{ep_0_order}
\begin{align}
-&\frac{\partial \sigma^0_{ij}}{\partial y_j} = 0 & \text{in } &Y_{f}, \label{ep_0_divstress} \\  
 &\sigma^0_{ij} = -p^0 \, \delta_{ij} + 2 \, \nu \, (e_{ijx}(\vc{v^0}) + e_{ijy}(\vc{v^1})) \label{sigma_0} \\
 &\frac{\partial v^0_j}{\partial x_j} + \frac{\partial v^1_j}{\partial y_j} = 0 & \text{in } &Y_{f}, \label{ep_0_order_incomp}\\
 &e_{ijx}(\vc{v^0}) + e_{ijy}(\vc{v^1}) = 0 & \text{in } &T, \label{ep_0_order_strain}\\
 &\frac{\partial B^0_j}{\partial y_j} = 0, \quad \ep_{ijk}\frac{\partial H^0_k}{\partial y_j} = 0 \text{ where } B^0_i = \mu H^0_i & \text{in } &Y, \label{cell_curl}
\end{align} 
\end{subequations}
with boundary conditions
\begin{gather}\label{ep_0_order_bc}
\begin{aligned}
\jump{ \vc{v^1} } &= \vc{0}, \quad \jump{ \vc{B^0} \cdot \vc{n} } = 0, \quad \jump{ \vc{n} \times \vc{H^0} } = \vc{0} \quad \text{ on $S$ }, \\
\quad &\vc{v^1}, \quad \vc{H^0} \text{ are } Y-\text{periodic}.
\end{aligned}
\end{gather}
Here $Y_f$ and $T$ denote the fluid, respectively the particle part of $Y$; and $S$ denotes the surface of $T$. At order of $\ep^2$ and $\ep^3$ we obtain from \eqref{bal_forces5} the balance of forces and torques for the particle $T$ respectively, 
\begin{gather}
\begin{aligned}
0 &= \int_S \sigma^0 \vc{n} \, ds + \alpha \int_S \jump{\tau^0 \vc{n}}) \, ds - \alpha \, \int_{T} \vc{B^0} \times {\rm curl_y{\vc{H}^0}} \, d\vc{y},\\ 
0 = \int_{S} \vc{y} &\times \sigma^0 \vc{n} \, ds + \alpha \, \int_S \vc{y} \times \jump{\tau^0 \vc{n}} \, ds - \alpha \, \int_{T} \vc{y }\times \left( \vc{B^0} \times {\rm curl_y(\vc{H}^0)} \right)  \, d\vc{y},
\end{aligned}
\end{gather}
where $\tau^0_{ij}$ is the Maxwell stress at order $\ep^0$ 

\begin{equation} \label{tau_0}
\tau_{ij}^0 = \mu \, H^0_i\, H^0_j - \frac{1}{2} \, \mu \, H^0_k\, H^0_k\, \delta_{ij},
\end{equation}

We remark that since from \eqref{cell_curl}  ${\rm curl_y (\vc{H^0})} = \vc{0}$ in $Y$, the balance of forces and torques simplify to the following,
\begin{equation} \label{cell_balance}
0 = \int_{S} \sigma^0 \vc{n} + \alpha \, \int_S \jump{\tau^0 \vc{n}} \, ds \text{ and } 0 = \int_{S} \vc{y} \times \sigma^0 \vc{n} \, ds + \alpha \, \int_S \vc{y} \times \jump{\tau^0 \vc{n}} \, ds. 
\end{equation}

\begin{remark}
At first order, the problem \eqref{ep_0_order}-\eqref{cell_balance} becomes one way coupled, as one could solve the Maxwell equations \eqref{cell_curl} independently. Once a solution is obtained, the Stokes problem \eqref{ep_0_divstress}-\eqref{ep_0_order_incomp}, can be solved with a known magnetic force added to the balance of forces and torques \eqref{cell_balance}.
\end{remark}

\section{Constitutive relations for Maxwell's equations}
\subsection{Study of the local problem}
Using the results from the two scale expansions, \eqref{cell_curl}, we can see that ${\rm curl_y (\vc{H^0})} = \vc{0}$ in $Y$ and thus there exists a function $\psi = \psi(\vc{x},\vc{y})$ with average $\widetilde{\psi}=0$ such that 

\begin{equation}\label{mag_pot}
H^0_i = -\frac{\partial \psi (\vc{x}, \vc{y})}{\partial y_i} + \widetilde{H}^0_i (\vc{x}) ,
\end{equation} 
where $\displaystyle \widetilde{\cdot} = \frac{1}{|Y|} \int_Y \cdot \, d\vc{y}$. Using the fact ${\rm div_y \vc{B^0}} = 0$ in $Y$, $B^0_i = \mu \, H^0_i$ and the boundary conditions \eqref{bc_five} we have,

\begin{gather} \label{psi_loc}
\begin{aligned}
-\frac{ \partial }{\partial y_i } \left ( \mu \, \left ( -\frac{ \partial \psi }{\partial y_i} + \widetilde{H}^0_i \right ) \right ) &= 0 \quad \text{ in $Y$},\\
\jump{ \mu \, \left ( -\frac{ \partial \psi }{\partial y_i} + \widetilde{H}^0_i \right ) \, n_i } &= 0 \quad \text{ on $S$ },\\
\psi \text{ is } Y-\text{periodic}, \quad &\widetilde{\psi} = 0.& 
\end{aligned} 
\end{gather}

Introducing the space of periodic functions, with zero average 

\[ \mathcal{W} = \left \{ w \in H^1_{per} (Y) \mid \, \widetilde{w} = 0 \right\}, \]
then the variational formulation of \eqref{psi_loc} is

\begin{gather}\label{var_loc_psi}
\begin{aligned}
		&\text{Find $\psi$ } \in \mathcal{W} \text{ such that} \\
		&\int_{Y} \mu \, \frac{ \partial \psi }{\partial y_i} \, \frac{ \partial v}{\partial y_i } \, d\vc{y}  = \widetilde{H}^0_i(\vc{x}) \, \int_{Y} \mu \,  \frac{ \partial v}{\partial y_i } \, d\vc{y} \text{ for any } v \in \mathcal{W}. 	
\end{aligned}
\end{gather}
Since we have imposed that $\psi$ has zero average over the unit cell $Y$, the solution to \eqref{var_loc_psi} can be determined uniquely by a simple application of the Lax-Milgram lemma.

Let $\phi^k$ be the unique solution of 
\begin{gather}\label{Maxwell_local}
\begin{aligned}
		&\text{Find $\phi^k$ } \in \mathcal{W} \text{ such that} \\
		&\int_{Y} \mu \, \frac{ \partial \phi^k }{\partial y_i} \, \frac{ \partial v}{\partial y_i } \, d\vc{y}  = \int_{Y} \mu \,  \frac{ \partial v}{\partial y_k } \, d\vc{y}\text{ for any } v \in \mathcal{W}. 	
\end{aligned}
\end{gather}

\begin{comment}
\begin{figure}[h]
    \centering
    \begin{minipage}{0.4\textwidth}
        \includegraphics[width=\textwidth]{phi1.png} % first figure itself
        \caption{Plot of $\phi^1$}
    \end{minipage}\hspace{1.0cm}
    \begin{minipage}{0.4\textwidth}
        \includegraphics[width=\textwidth]{phi2.png} % second figure itself
        \caption{Plot of $\phi^2$}
    \end{minipage}
\caption{Plot of the local solution $\phi^k$ in \eqref{Maxwell_local} for iron particles of $19\%$ volume fraction using $50 \times 50$ $P1$ elements to mesh the square and $100$ $P1$ elements to mesh the circle. The effective magnetic permeability of the particle was set to $2\times10^5$ while the effective magnetic permeability of the fluid was set to $1.0$. Images were obtained using {\tt FreeFem++}~\cite{FH12}.}
\end{figure}
\end{comment}

By virtue of linearity of \eqref{var_loc_psi} we can write,
 
$$\psi(\vc{x},\vc{y}) = \phi^k (\vc{y}) \, \widetilde{H}^0_k (\vc{x}) + C(\vc{x}).$$ 

In principle, once $\widetilde{H}^{0}_{k}$ is known, we can determine $\psi$ up to an additive function of $\vc{x}.$ Hence, combining \eqref{mag_pot} and the above relationship between $\psi$ and $\phi^k$ we obtain the following constitutive law between the magnetic induction and the magnetic field,

\begin{equation} \label{Max_const}
	\widetilde{B}^0_i = \mu^H_{ik} \, \widetilde{H}^0_k, \text{ where } \mu^H_{ik} = \int_Y \mu \, \left ( -\frac{ \partial \phi^k }{\partial y_i} + \delta_{ik} \right ) \, d\vc{y}.
\end{equation}

One can show (see \cite{SP80}) that the homogenized magnetic permeability tensor is symmetric, $\mu^H_{ik} = \mu^H_{ki}$. Moreover, if we denote by $A_{i\ell}(\vc{y}) = \left ( -\frac{ \partial \phi^\ell(\vc{y}) }{\partial y_i} + \delta_{i\ell} \right )$ one can see from \eqref{mag_pot} that $H^0_i =  A_{i\ell} \widetilde{H}^0_\ell$ and thus the Maxwell stress \eqref{tau_0} takes the following form,

\begin{align}
\tau_{ij}^0 = \mu \, A_{i\ell} \, A_{jm} \, \widetilde{H}^0_\ell \, \widetilde{H}^0_m - \frac{1}{2} \, \mu \, A_{mk} \, A_{\ell k} \, \delta_{ij} \, \widetilde{H}^0_m \, \widetilde{H}^0_\ell = \mu \, A_{ij}^{m\ell} \, \widetilde{H}^0_m \, \widetilde{H}^0_\ell. \nonumber
\end{align}

Here $A_{ij}^{m\ell} =  \frac{1}{2} \, \left( A_{i\ell} \, A_{jm} + A_{j\ell} \, A_{im} - A_{mk} \, A_{\ell k} \, \delta_{ij} \right)$ and has the following symmetry, $A_{ij}^{m\ell} = A_{ji}^{m\ell}=A_{ij}^{\ell m}$. Recall that the $\dv{\tau^\ep} = 0$ in $\Om_{2\,\ep}$ and $\dv{\tau^\ep} = -\vc{B^\ep} \times \curl{\vc{H^\ep}}$ in $\Om_{1\ep}$. From the two scale expansion, at order $\ep^{-1}$ from equation \eqref{mstress_div} we obtain, 

\begin{equation} \label{dv_tau_0}
\dv_y{\tau^0} = 0 \text{ in } Y.
\end{equation}

\section{Fluid velocity and pressure}
\subsection{Study of the local problems}

Problem \eqref{ep_0_order}-\eqref{ep_0_order_bc}, \eqref{cell_balance} is an elliptic problem in the variable $\vc{y} \in Y$ with forcing terms $\vc{v^0}(\vc{x})$ and $\widetilde{\vc{H}}^0(\vc{x})$ at the macroscale. We can decouple the contributions of $\vc{v^0}(\vc{x})$ and $\widetilde{H}^0(\vc{x})$ and split $\vc{v^1}$ and $p^0$ in two parts: a part that is driven by the bulk velocity, and a part that comes from the bulk magnetic field.

\begin{equation} \label{v1}
	v^1_k(\vc{x},\vc{y}) = \chi^{m\ell}_k(\vc{y}) \, e_{m\ell}(\vc{v^0}) + \xi^{m\ell}_k(\vc{y}) \, \widetilde{H}^0_m \, \widetilde{H}^0_\ell + A_k(\vc{x}),
\end{equation}
\begin{equation} \label{p0}
p^0(\vc{x},\vc{y}) = p^{m\ell}(\vc{y})\,e_{m\ell}(\vc{v^0}) + \pi^{m\ell}(\vc{y})\,\widetilde{H}^0_m \, \widetilde{H}^0_\ell + \bar{p}^0(\vc{x}),
\end{equation}
where $\displaystyle \int_{Y_f} p^{m\ell}(\vc{y}) \, d\vc{y}=0$ and $\displaystyle \int_{Y_f} \pi^{m\ell}(\vc{y}) \, d\vc{y}=0.$ 

Here, $\vc{\chi^{ml}}$ satisfies 
\begin{gather} \label{chi_P}
\begin{aligned}
-\frac{ \partial }{\partial y_j } \varepsilon_{ij}^{m\ell} &= 0 \quad \text{ in $Y_f$},\\
\varepsilon_{ij}^{m\ell} = -p^{m\ell} \delta_{ij} &+ 2 \, (C_{ijm\ell} + e_{ijy}(\vc{\chi}^{m\ell})) \\ 
-\frac{ \partial \chi^{m\ell}_i}{\partial y_i } &= 0 \quad \text{ in $Y_f$}, \\
\jump{ \vc{\chi}^{m\ell}} &= 0 \quad \text{ on $S$ },\\
C_{ijm\ell} + e_{ijy}(\vc{\chi}^{m\ell}) &= 0 \quad \text{ in $T$ },\\
\vc{\chi}^{m\ell} \text{ is } & Y-\text{periodic}, \quad \widetilde{\vc{\chi}}^{m\ell}=\vc{0} \text{ in } Y, 
\end{aligned} 
\end{gather}
together with the balance of forces and torques,

\begin{equation}\label{FTF}
\int_{S} \varepsilon_{ij}^{m\ell} n_j \, ds = 0 \text{ and } \int_{S} \ep_{ijk}\,y_j\,\varepsilon_{kp}^{m\ell} n_p \, ds = 0,
\end{equation}
where $\displaystyle C_{ijm\ell} = \frac{1}{2}(\delta_{im}\delta_{j\ell} + \delta_{i\ell}\delta_{jm}) - \frac{1}{n}\delta_{ij}\,\delta_{m\ell}$. Equation~\eqref{chi_P} is well known having been obtained by many authors \cite{LSP185, LSP285, LV94} among others. 

The variational formulation of \eqref{chi_P}-\eqref{FTF} is:

\begin{gather}\label{var_chi}
\begin{aligned}
		&\text{Find $\vc{\chi^{m\ell}}$ } \in \mathcal{U} \text{ such that} \\
		&\int_{Y_f} 2 \, e_{ij}(\vc{\chi^{m\ell}}) \, e_{ij}(\vc{\phi}-\vc{\chi^{m\ell}}) \, d\vc{y} = 0, \text{ for all } \vc{\phi} \in \mathcal{U}, 	
\end{aligned}
\end{gather}
where $\mathcal{U}$ is the closed, convex, non-empty subset of $H_{per}^1(Y)^n$ defined by

\begin{align}
\mathcal{U} = \left \{  \vc{u} \in H_{per}^1(Y)^n \mid \dv{\vc{u}}=0 \text{ in } Y_f, e_{ij}(\vc{u})=-C_{ijm\ell} \text{ in } T, \jump{\vc{u}}=\vc{0} \text{ on } S, \right. \nonumber \\
\left. \widetilde{\vc{u}} = \vc{0} \text{ in } Y \right \}. \nonumber
\end{align}

\begin{remark}
We remark that if we define $B^{ij}_k = \frac{1}{2}(y_i\, \delta_{jk} + y_j\,\delta_{ik}) - \frac{1}{n}y_k\,\delta_{ij}$, then it immediately follows that $e_{ij}(\vc{B}^{m\ell}) = C_{ijm\ell}$. 
\end{remark}
Existence and uniqueness of a solution follows from classical theory of variational inequalities \cite{KinSta}. In similar fashion we can derive the local problem for $\vc{\xi^{ml}}$, 

\begin{gather}\label{xi_P}
\begin{aligned}
-\frac{ \partial }{\partial y_j } \Sigma_{ij}^{m\ell} &= 0 \quad \text{ in $Y_f$},\\
\Sigma_{ij}^{m\ell} = -\pi^{m\ell} \delta_{ij} &+ 2 \, e_{ijy}(\vc{\xi}^{m\ell})  \\ 
-\frac{ \partial \xi^{m\ell}_i}{\partial y_i } &= 0 \quad \text{ in $Y_f$}, \\
\jump{ \vc{\xi}^{m\ell}} &= 0 \quad \text{ on $S$ },\\
e_{ijy}(\vc{\xi}^{m\ell}) &= 0 \quad \text{ in $T$ },\\
\vc{\xi}^{m\ell} \text{ is } Y-&\text{periodic}, \quad \widetilde{\vc{\xi}}^{m\ell}=0.
\end{aligned} 
\end{gather}

Using \eqref{dv_tau_0} the balance of forces reduces to,
 
\begin{equation}
\int_{S} \Sigma_{ij}^{m\ell} \, n_j \, ds = 0,
\end{equation} 
together with the balance of torques

\begin{equation} \label{FTM} 
\int_{S} \ep_{ijk}\,y_j\,\left ( \Sigma_{kp}^{m\ell} + \alpha \jump{\mu A_{kp}^{m\ell}} \right ) \, n_p \, ds = 0.
\end{equation}

We can formulate \eqref{xi_P}--\eqref{FTM} variationally as,
 
\begin{gather}\label{var_xi}
\begin{aligned}
		&\text{Find $\vc{\xi^{m\ell}}$ } \in \mathcal{V} \text{ such that} \\
		&\int_{Y_f} 2 \, e_{ijy}(\vc{\xi^{m\ell}}) \, e_{ijy}(\vc{\phi}) \, d\vc{y} +  \int_{Y} A_{ij}^{m\ell} \, e_{ijy}(\vc{\phi}) \, d\vc{y} = 0, \text{ for all } \vc{\phi} \in \mathcal{V}, 	
\end{aligned}
\end{gather}
where 

\[
%\begin{split}
\mathcal{V}= \left \{ \vc{v} \in H_{per}^1(Y)^n \mid \dv{\vc{v}}=0 \text{ in } Y_f, e_{y}(\vc{u})=0 \text{ in } T, \jump{\vc{v}}=\vc{0} \text{ on } S, \widetilde{\vc{v}} = \vc{0} \text{ in } Y \right \},
%\end{split}
\]
is a closed subspace of $H_{per}^1(Y)^n$. Existence and uniqueness follows from an application of the Lax-Milgram lemma. These equations indicate the contribution of the magnetic field and the solution $\vc{\xi^{m\ell}}$ depends, through the balance of forces and torques on the solution of the local problem~\eqref{Maxwell_local} and the effective magnetic permeability of the composite.

\begin{remark}
We remark that the only driving force that makes the solution $\vc{\xi^{m\ell}}$ non trivial in \eqref{var_xi} is the rotation induced by the magnetic field through the fourth order tensor $A_{ij}^{m\ell}$.
\end{remark}

\section{Effective balance equations}

The two-scale expansion at the $\ep^0$ order yields the following problems:
\begin{subequations}
\begin{align}
-\dv_x{\sigma^0} - \dv_y{\sigma^1} &= \vc{0} \quad &\text{ in $Y_f$},\label{hom_fluid_stress} \\ 
\dv_x{\vc{v^1}} + \dv_y{\vc{v^2}} &= 0 \quad &\text{ in $Y_f$}, \\
\dv_x{\vc{B^0}} + \dv_y{\vc{B^1}} &= 0 \quad &\text{ in $Y$}, \label{hom_induct} \\
\curl_x{\vc{H^0}} + \curl_y{\vc{H^1}} &= 0 \quad &\text{ in $Y_f$}, \label{hom_mag_field1} \\
\curl_x{\vc{H^0}} + \curl_y{\vc{H^1}} &= R_m \, \vc{v^0} \times \vc{B^0} \quad &\text{ in $T$} \label{hom_mag_field2},
\end{align} 
\end{subequations}
with boundary conditions 

\begin{gather}
\begin{aligned}
\jump{ \vc{v^2} } &= \vc{0}, \quad \jump{ \vc{B^1} \cdot \vc{n} } = 0 \quad \jump{ \vc{n} \times \vc{H^1} } = \vc{0} \quad \text{ on $S$ }, \\
&\vc{v^2}, \quad \vc{H^1} \text{ are } Y-\text{periodic}.
\end{aligned}
\end{gather}

In each period, we consider a Taylor expansion, around the center of mass of the rigid particle, both of the viscous stress and the Maxwell stress of the form (see \cite{LV94}),
\[
\sigma^\ep (\vc{x}) = \sigma^0 (\vc{x^\ell_c}, \vc{y}) + \frac{\partial \sigma^0 (\vc{x^\ell_c}, \vc{y})}{\partial x_\alpha}(x_\alpha - x^\ell_{c, \alpha}) + \ep \, \sigma^1 (\vc{x^\ell_c}, \vc{y}) + \ep \, \frac{\partial \sigma^1 (\vc{x^\ell_c}, \vc{y})}{\partial x_\alpha}(x_\alpha - x^\ell_{c, \alpha}) + \cdots
\]

\[
\tau^\ep (\vc{x}) = \tau^0 (\vc{x^\ell_c}, \vc{y}) + \frac{\partial \tau^0 (\vc{x^\ell_c}, \vc{y})}{\partial x_\alpha}(x_\alpha - x^\ell_{c, \alpha}) + \ep \, \tau^1 (\vc{x^\ell_c}, \vc{y}) + \ep \, \frac{\partial \tau^1 (\vc{x^\ell_c}, \vc{y})}{\partial x_\alpha}(x_\alpha - x^\ell_{c, \alpha}) + \cdots
\]
where the expansion of the Maxwell stress occurs both inside the rigid particle and the fluid. Using this method we can expand the balance of forces, \eqref{bal_forces5}, and obtain at order $\ep^3$,
\begin{gather}\label{stress_exp}
\begin{aligned} 
0 &= \int_{S} \left ( \frac{\partial \sigma^0_{ij} }{ \partial x_k } \, y_k + \sigma^1_{ij}\right ) \, n_j \, ds + \alpha \, \int_S \jump{\left(\frac{\partial \tau^0_{ij} }{ \partial x_k } \, y_k + \tau^1_{ij} \right) \, n_j} \, ds \\
&- \alpha \, \int_{T} (\vc{B^0} \times (\curl_x{\vc{H^0}} + \curl_y{\vc{H^1}}))_i \, d\vc{y}. 
\end{aligned}
\end{gather}

Integrate \eqref{hom_fluid_stress} over $Y_f$ and add to \eqref{stress_exp} obtain the following,
\begin{gather}\label{hom_stress_exp}
\begin{aligned} 
0 &= \int_{Y_f} \frac{\partial \sigma^0_{ij} }{ \partial x_j } \, d\vc{y} + \int_{S} \frac{\partial \sigma^0_{ij} }{ \partial x_k } \, y_k n_j \, ds + \alpha \, \int_S \jump{\right(\frac{\partial \tau^0_{ij} }{ \partial x_k } \, y_k + \tau^1_{ij} \left) \, n_j} \, ds \\
&- \alpha \, \int_{T} (\vc{B^0} \times (\curl_x{\vc{H^0}} + \curl_y{\vc{H^1}}))_i \, d\vc{y}.
\end{aligned}
\end{gather}

At order $\ep^0$ we obtain, $\dv_x{\tau^0} + \dv_y{\tau^1} = 0$  in $Y_f$ and $\dv_x{\tau^0} + \dv_y{\tau^1} = -\vc{B^0} \times (\curl_x{\vc{H^0}} + \curl_y{\vc{H^1}})$ in $T$. Combining the aforementioned results and the divergence theorem we can rewrite \eqref{hom_stress_exp} the following way,
\begin{align} \label{macro_equation1}
0 &= \int_{Y_f} \frac{\partial \sigma^0_{ij} }{ \partial x_j } \, d\vc{y} + \int_{S} \frac{\partial \sigma^0_{ik} }{ \partial x_j } \, y_j n_k \, ds + \alpha \, \int_S \jump{\frac{\partial \tau^0_{ik} }{ \partial x_j } \, y_j \, n_k} \, ds + \alpha \, \int_Y  \frac{\partial \tau^0_{ij} }{ \partial x_j } \, d\vc{y}.
\end{align}

Using the decomposition of $\vc{v^1}$ and $p^0$ in \eqref{v1} and \eqref{p0} we can re-write $\sigma^0_{ij}$ and $\tau^0_{ij}$, 
\[
\sigma_{ij}^0 = -\bar{p}^0 \, \delta_{ij} + \varepsilon_{ij}^{m\ell}e_{mlx}(\vc{v^0}) + \Sigma_{ij}^{m\ell}\widetilde{H}^0_m \, \widetilde{H}^0_\ell, \quad
\tau_{ij}^0 = \mu\,A_{ij}^{m\ell}\widetilde{H}^0_m \, \widetilde{H}^0_\ell.
\]

Moreover, equations \eqref{sigma_0}, \eqref{tau_0}, \eqref{chi_P} and \eqref{xi_P} allow us to retain the only symmetric part of \eqref{macro_equation1}. Hence the homogenized fluid equations \eqref{macro_equation1} become,
\begin{multline} \label{macro_equation2}
0 = \frac{\partial}{\partial x_j}  \Bigl( -\bar{p}^0 \delta_{ij} + \Bigl \{ \int_{Y_f} 2 e_{ijy}(\vc{B^{m\ell}} + \vc{\chi^{m\ell}}) d\vc{y} + \int_S \varepsilon_{pk}^{m\ell}  B^{ij}_p \, n_k ds \Bigr \} e_{m\ell x}(\vc{v^0}) \\
+ \left \{ \int_{Y_f} 2 e_{ijy}(\vc{\xi^{m\ell}}) d\vc{y} + \int_S \Sigma_{pk}^{m\ell} B^{ij}_p n_k ds \right. \\ 
\left. + \alpha \int_Y \mu A_{ij}^{m\ell} d\vc{y} + \alpha \int_S \jump{\mu A_{pk}^{m\ell}} B^{ij}_p n_k \right \}   \widetilde{H}^0_m \widetilde{H}^0_\ell \Bigl ). 
\end{multline}
Furthermore, using \eqref{ep_0_order_incomp}--\eqref{ep_0_order_strain} and the divergence theorem we can obtain the incompressibility condition, $\dv_x{\vc{v^0}} = 0$. 

Denote by

\[
\nu^H_{ijm\ell} =  \left \{ \int_{Y_f} 2 \,e_{ijy}(\vc{B^{m\ell}} + \vc{\chi^{m\ell}}) \, d\vc{y} + \int_S \varepsilon_{pk}^{m\ell} \, B^{ij}_p \, n_k \, ds \right \}
\]
and

\[
\beta^H_{ijm\ell} = \left \{ \int_{Y_f} 2 e_{ijy}(\vc{\xi^{m\ell}}) d\vc{y} + \int_S \Sigma_{pk}^{m\ell} B^{ij}_p n_k ds + \alpha \int_Y \mu A_{ij}^{m\ell} d\vc{y} + \alpha \int_S \jump{\mu A_{pk}^{m\ell}} B^{ij}_p n_k \right \}
\]
then the homogenized equation \eqref{macro_equation2} becomes

\[
0 = \frac{\partial}{\partial x_j} \left ( -\bar{p}^0 \, \delta_{ij} + \nu_{ijm\ell} \, e_{m\ell x}(\vc{v^0}) + \beta_{ijm\ell} \, \widetilde{H}^0_m \, \widetilde{H}^0_\ell \right).
\]

Using local problem \eqref{chi_P} we can re-write the $\nu_{ijm\ell}$ the following way,
\begin{equation} \label{nu}
\nu^H_{ijm\ell} = \int_{Y_f} 2 \, e(\vc{B^{ml}} + \vc{\chi^{ml}}) : e(\vc{B^{ij}} + \vc{\chi^{ij}}) \, d\vc{y}.
\end{equation}
which is a well known formula derived in~\cite{LSP185}, \cite{LSP285}, \cite{SP80} as well as its generalizations derived in~\cite{LV94}, \cite{NV15}. In a similar fashion, using local problem \eqref{xi_P} and the kinematic condition in \eqref{chi_P} we can re-write $\beta^H_{ijm\ell}$ as follows,
\begin{gather} \label{beta}
\begin{aligned}
\beta^H_{ijm\ell} =  \int_{Y_f} 2 e(\vc{\xi^{ml}}):e(\vc{B^{ij}} + \vc{\chi^{ij}}) d\vc{y} &+ \alpha \int_{Y_f} \mu A^{m\ell}:e(\vc{B^{ij}} + \vc{\chi^{ij}}) d\vc{y}\\ &+ \alpha \int_Y \mu A_{ij}^{m\ell} d\vc{y}.
\end{aligned}
\end{gather}

It is now clear that $\nu^H_{ijm\ell}$ possesses the following symmetry, $\nu^H_{ijm\ell} = \nu^H_{jim\ell} = \nu^H_{m \ell i j}$. While for $\beta^H_{ijm\ell}$, we have $\beta^H_{ijm\ell} = \beta^H_{jim\ell}=\beta^H_{i j \ell m}.$ 

To obtain the homogenized Maxwell equations, average \eqref{hom_induct}, \eqref{hom_mag_field1}, and \eqref{hom_mag_field2} over $Y$, $Y_f$, and $T$ respectively and use equation \eqref{Max_const} to obtain,

\begin{gather}
\begin{aligned} \label{hom_maxwell}
\frac{\partial \, (\mu^H_{ik} \, \widetilde{H}^0_k)}{\partial x_j} = 0, \quad \ep_{ijk} \frac{\partial \widetilde{H}^0_k}{\partial x_j} = R_m \, \ep_{ijk} \, v^0_j \, \mu^{HS}_{kp} \, \widetilde{H}^0_p \quad \text{ in } \Om, \nonumber
\end{aligned}
\end{gather}
where 
\begin{equation} \label{mu:s} 
\mu^{HS}_{ik} = \int_T \mu \, \left ( -\frac{ \partial \phi^k }{\partial y_i} + \delta_{ik} \right ) \, d\vc{y} 
\end{equation} 
with boundary conditions,
\begin{gather*}
\begin{aligned} 
\widetilde{H}^0_i = b_i, \quad v^0_i =0 \quad \text{ on } \partial \Om. 
\end{aligned}
\end{gather*}

The effective coefficients are computed as the angular averaging of the tensors $\nu^H_{ijm\ell}$ and $\beta^H_{ijm\ell}$. This is done by introducing the projection on hydrostatic fields, $P_b$, and the projection on shear fields $P_s$ (see~\cite{MV10}, \cite{Mil}, \cite{MK88}). The components of the projections in three dimensional space are given by: 
\[
(P_b)_{ijk\ell} = \frac{1}{n} \, \delta_{ij} \, \delta_{k\ell}, \quad (P_s)_{ijk\ell} = \frac{1}{2} \, (\delta_{ik} \, \delta_{j\ell} + \delta_{i\ell} \, \delta_{jk}) - \frac{1}{n} \, \delta_{ij} \, \delta_{k\ell}.
\]

Let us fix the following notation:
\begin{align}
\nu_b=tr(P_b \, \nu^H) &= \frac{1}{n} \, \nu^H_{ppqq}, \quad \nu_s = tr(P_s \, \nu^H) = \left( \nu^H_{pqpq} -\frac{1}{n} \, \nu^H_{ppqq} \right), \nonumber \\ 
\beta^H_b = tr(P_b \, \beta^H) &= \frac{1}{n} \, \beta^H_{ppqq}, \quad \beta_s = tr(P_s \, \beta^H) = \left( \beta^H_{pqpq} -\frac{1}{n} \beta^H_{ppqq}\right).\nonumber
\end{align}
we can re-write the homogenized coefficients $\nu^H_{ijm\ell}$ and $\beta^H_{ijm\ell}$ as follows,
\[
\nu^H_{ijm\ell} = \frac{1}{n} \, (\nu_b - \nu_s) \delta_{ij}\delta_{m\ell} +  \frac{1}{2} \, \nu_s \, (\delta_{ik}\delta_{j\ell} + \delta_{i\ell}\delta_{jk}),
\]
\[
\beta^H_{ijm\ell} =  \frac{1}{n} \, (\beta_b - \nu_s) \delta_{ij}\delta_{m\ell} +  \frac{1}{2} \, \beta_s \, (\delta_{ik}\delta_{j\ell} + \delta_{i\ell}\delta_{jk}).
\]

Gathering all the equations we have that the homogenized equations governing the magnetorheological fluid form the following coupled system between the Stokes equations and the quasistatic Maxwell equations,
%\begin{gather} \label{fin_eq}
%\begin{aligned} 
%& \frac{\partial}{\partial x_j} \left ( \sigma^H_{ij} + \tau^H_{ij} \right ) = 0 \text{ in } \Om, \\
%& \sigma^H_{ij} + \tau^H_{ij} = -\bar{p}^0 \, \delta_{ij} + \nu_s \, e_{ij}(\vc{v^0}) + \frac{1}{n} \, (\beta_b - \beta_s) \, \delta_{ij} \left|\widetilde{\vc{H}}^0\right|^2 + \beta_s \, \widetilde{H}^0_i \, \widetilde{H}^0_j, \\
%& \frac{\partial{v_i^0}}{\partial x_i} =0 \text{ in } \Om, \\
%& \frac{\partial (\mu^H_{jk} \, \widetilde{H}^0_k)}{\partial x_j} = 0 \text{ in } \Om, \\ 
%& \ep_{ijk} \frac{\partial \widetilde{H}^0_k}{\partial x_j} = R_m \, \ep_{ijk} \, v^0_j \, \mu^{HS}_{kp} \, \widetilde{H}^0_p \text{ in } \Om, \\
%& v_i^0 = 0, \text{ on } \partial \Om, \\
%& \widetilde{H}_i^0 = b_i \text{ on } \partial \Om.
%\end{aligned}
%\end{gather}
%or written in compact vector notation we have, 
\begin{gather} \label{eq:vecform}
\begin{aligned} 
& \dv{ \left ( \sigma^H + \tau^H \right ) } = \vc{0} \text{ in } \Om, \\
& \sigma^H + \tau^H = -\bar{p}^0 \, I + \nu_s \, e(\vc{v^0}) + \frac{1}{n} \, (\beta_b - \beta_s) \, \left|\widetilde{\vc{H}}^0\right|^2 \, I + \beta_s \, \widetilde{\vc{H}}^0 \otimes \widetilde{\vc{H}}^0, \\
& \dv{\vc{v}^0} = 0 \text{ in } \Om, \\
& \dv{\left ( \mu^H \, \widetilde{\vc{H}}^0 \right )} = 0 \text{ in } \Om, \\ 
& \curl{\widetilde{\vc{H}}^0} = R_m \, \vc{v}^0 \, \times \mu^{HS} \, \widetilde{\vc{H}}^0 \text{ in } \Om, \\
& \vc{v}^0 = \vc{0}, \text{ on } \partial \Om, \\
& \widetilde{\vc{H}}^0 = \vc{b} \text{ on } \partial \Om.
\end{aligned}
\end{gather}

Equation~\eqref{eq:vecform} generalizes the quasistatic set of equations introduced in~\cite{NR64}, \cite{Ro14} in two ways: first by providing exact formulas for the effective coefficients which consist of the homogenized viscosity, $\nu^H$, and three homogenized magnetic permeabilities, $\mu^H$, $\mu^{HS}$, and $\beta^H$, which all depend on the geometry of the suspension, the volume fraction, the magnetic permeability $\mu$, the Alfven number $\alpha$, and the particles distribution. Second, by coupling the fluid velocity field with the magnetic field through Ohm's law.

\section{Numerical results for a suspension of iron particles in a viscous non-conducting fluid}

The goal of this section is to carry out calculations, using the finite element method, of the effective viscosity $\nu^H$ and the effective magnetic coefficients $\beta^H$, $\mu^H$, and $\mu^{HS}$ that describe the behavior of a magnetorheological fluid in \eqref{eq:vecform}. To achieve this we need to compute the solutions $\phi^k$, $\chi^{m\ell}$, and $\xi^{m\ell}$ of the local problems \eqref{Maxwell_local}, \eqref{chi_P}, and \eqref{xi_P} respectively. 

Unlike regular suspensions for which the effective properties are dependent only on fluid viscosity, particle geometry, and volume fraction, for magnetorheological fluids of significance is also the particles' distribution. The magnetic field polarizes the particles which align in the field direction to form {\it chains} and {\it columns} and that contributes significantly to the increase of the yield stress~\cite{Bossis}, \cite{tao01}, \cite{Win49}. 

The choice of the periodic unit cell, as well as the geometry and distribution of particles, can lead to different {\it chain structures} and hence different effective properties. We are considering here a uniform distribution of particles as well as a {\it chain} distribution; to achieve this we change the aspect ratio of the unit cell from a $1 \times 1$ square to a $2 \times 1/2$ rectangle. 

In all computations we have used a regular, symmetric, triangular mesh. For the standard $1 \times 1$ periodic unit cell we used $20 \times 20$ $P1$ elements to mesh the square and $100$ $P1$ elements to mesh the circular particle inside the square, while for the elongated $2 \times \frac{1}{2}$ periodic unit cell, we used $100 \times 20$ $P1$ elements for the rectangle and $100$ $P1$ elements for the circular particle. Other geometries where the particle is an ellipse, for instance, present an interesting case on their own. Ellipses have a priori a preferred direction i.e. they are anisotropic and as result the effective coefficients will be anisotropic. Hence, when one discusses {\it chain structures} of ellipses {\it angle orientation} must be taken into account. We do not explore such cases in the current work and we consider only the case of circular particles.

We remark that in the two dimensional setting, the tensors entries $C_{ijmm} = 0$ and $\vc{B}^{mm} = \vc{0}$. As a consequence, of the linearity of the local problem \eqref{chi_P}, we have $\vc{\chi}^{mm} = \vc{0}$. Hence, $\nu_{mmii} =0$ which implies that $\nu_b=0$. Using a similar argument, we can similarly show that $\beta_{mmii} = 0$ which implies that $\beta_b=0$. 

The relative magnetic permeability of the iron particle ($99.95\%$ pure) was fixed through out to be $2 \times 10^5$ while that of the fluid was set to $1$ \cite{CO58}. All the calculations were carried out using the software {\tt FreeFem++}~\cite{FH12}. 

\begin{figure}[!htb]
\label{fig:obo}
\centering
\begin{tabular}{ccc}
\subf{\includegraphics[width=45mm]{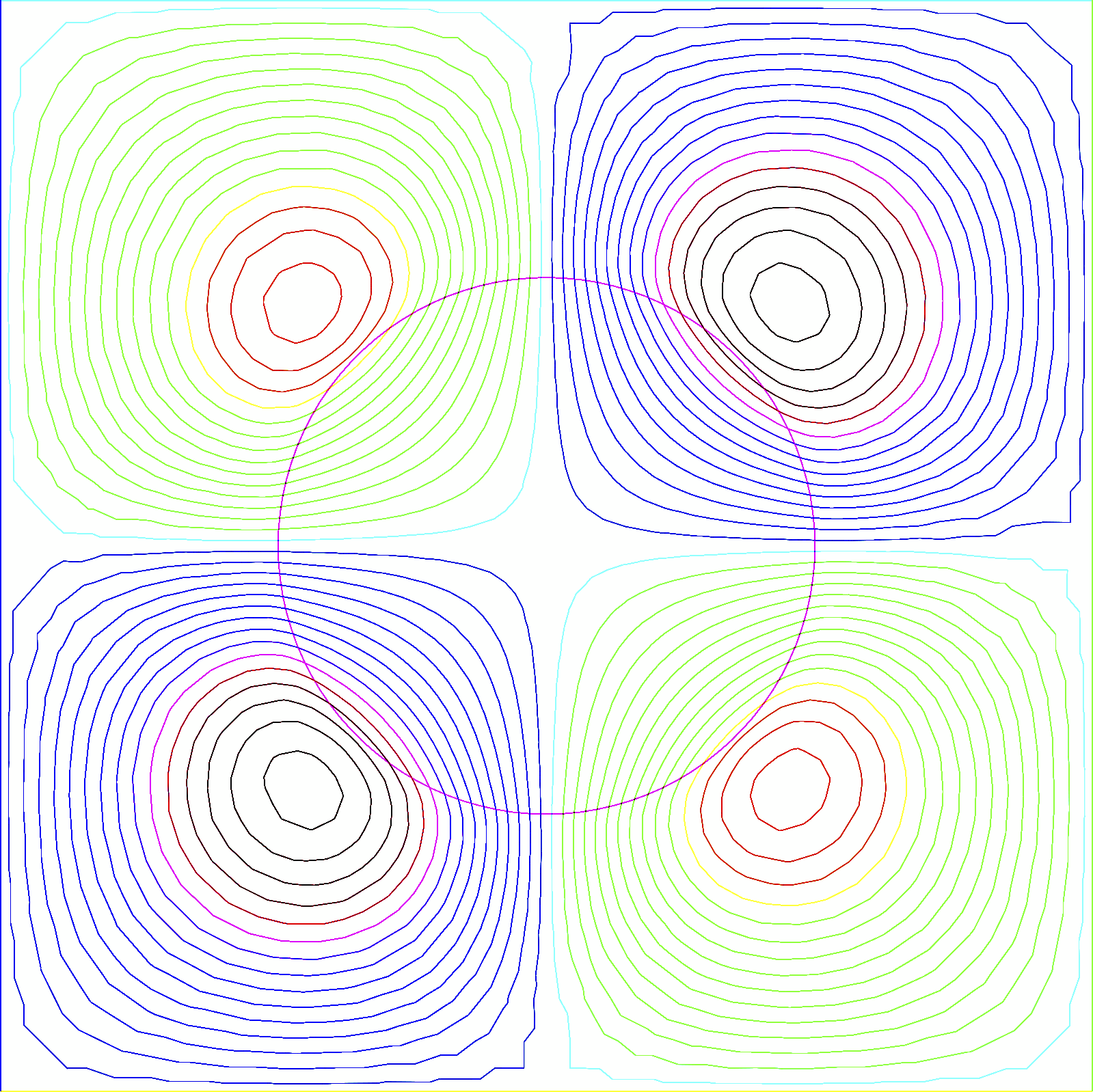}}
     {$\chi^{11}$}
&
\subf{\includegraphics[width=45mm]{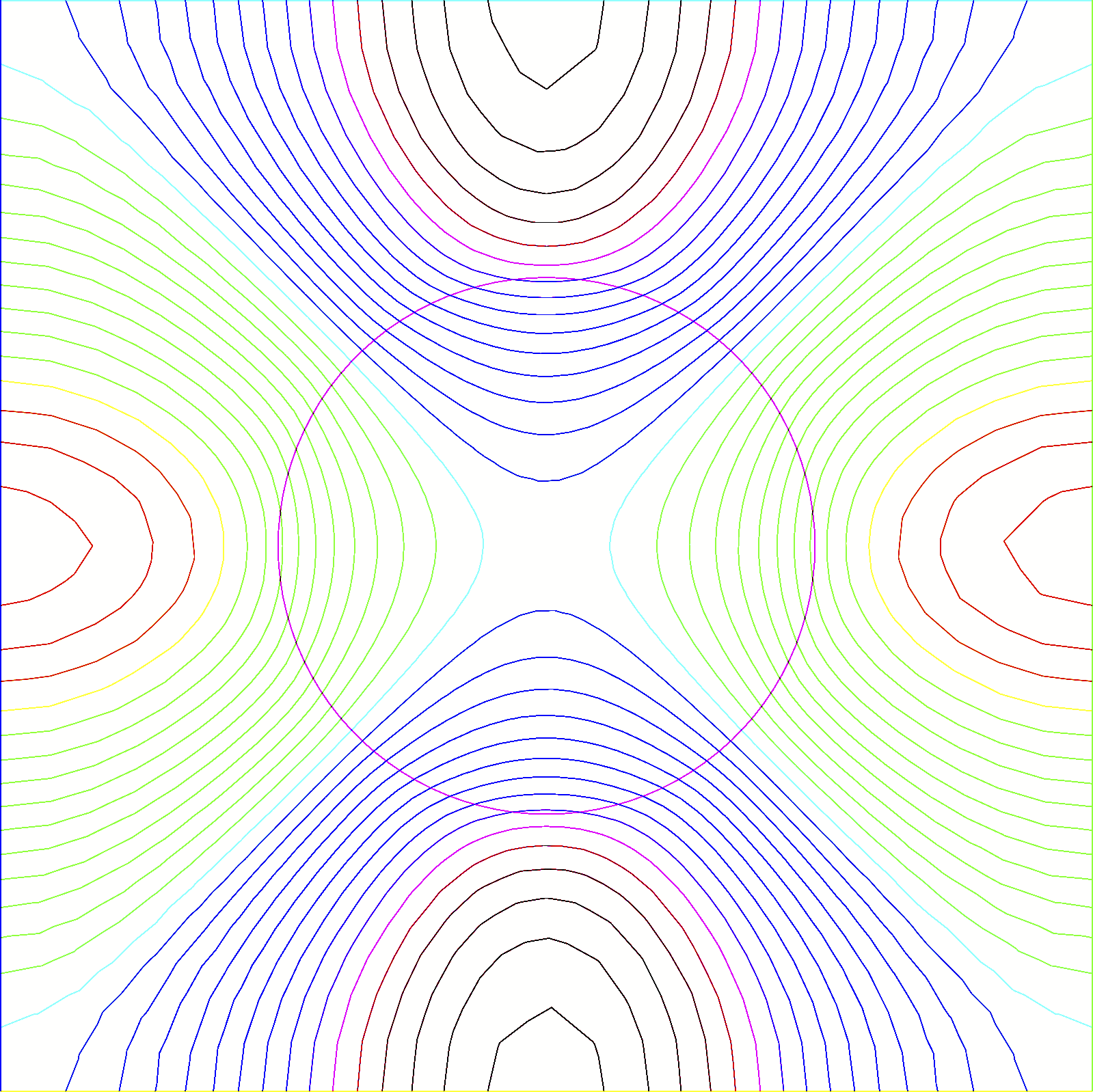}}
     {$\chi^{12}$}
&
\subf{\includegraphics[width=45mm]{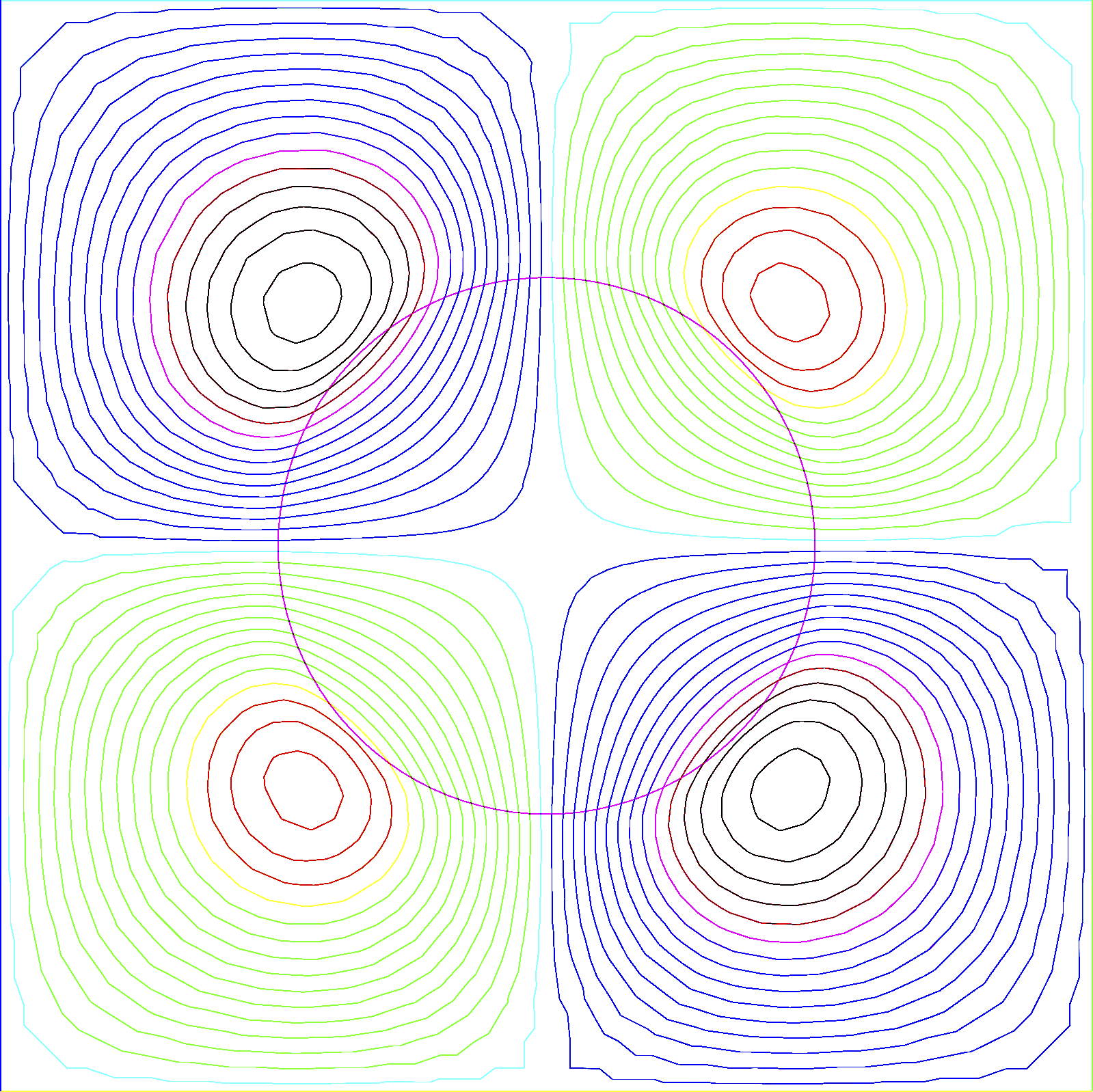}}
     {$\chi^{22}$}
\\
\subf{\includegraphics[width=45mm]{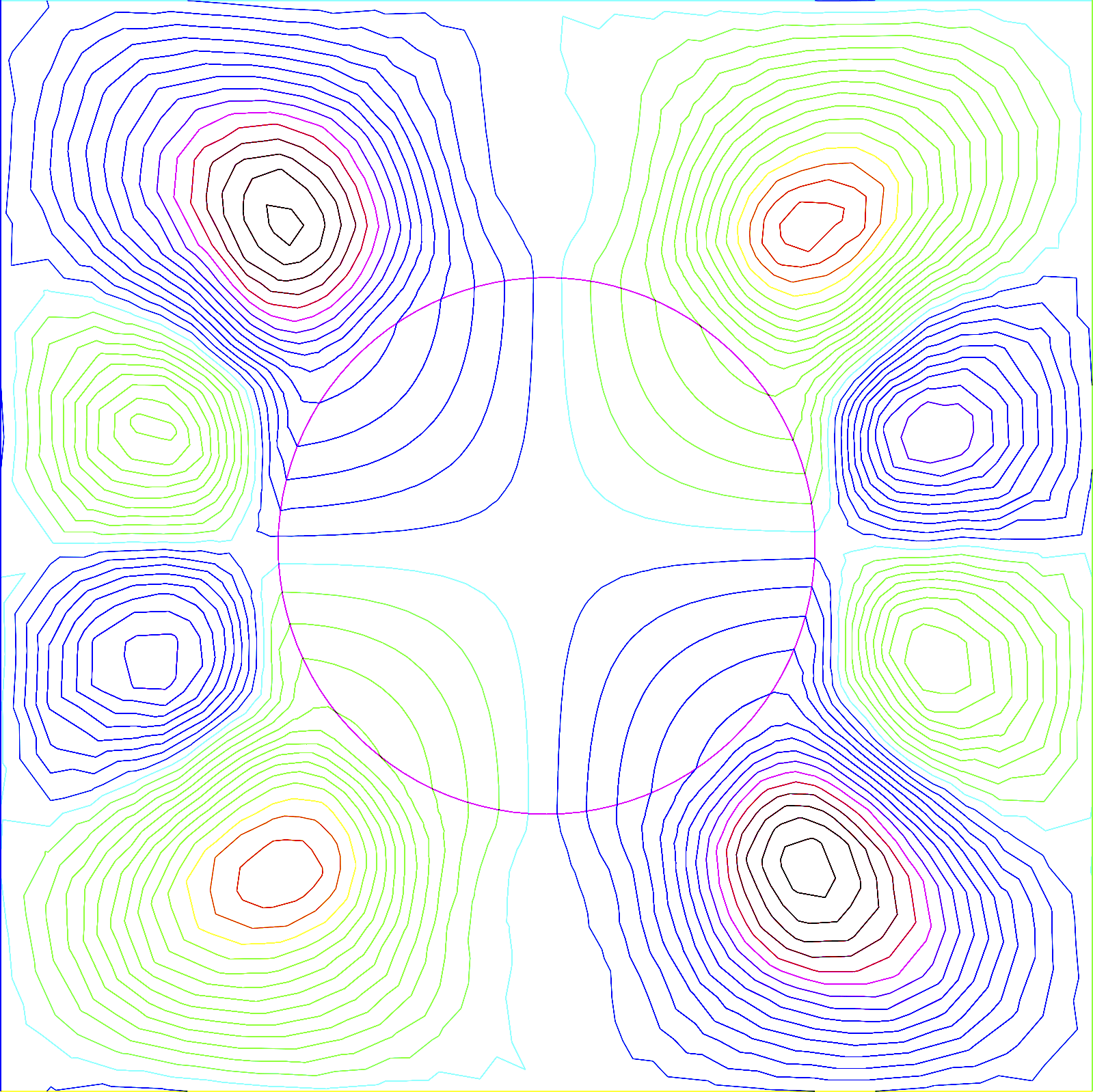}}
     {$\xi^{11}$}
&
\subf{\includegraphics[width=45mm]{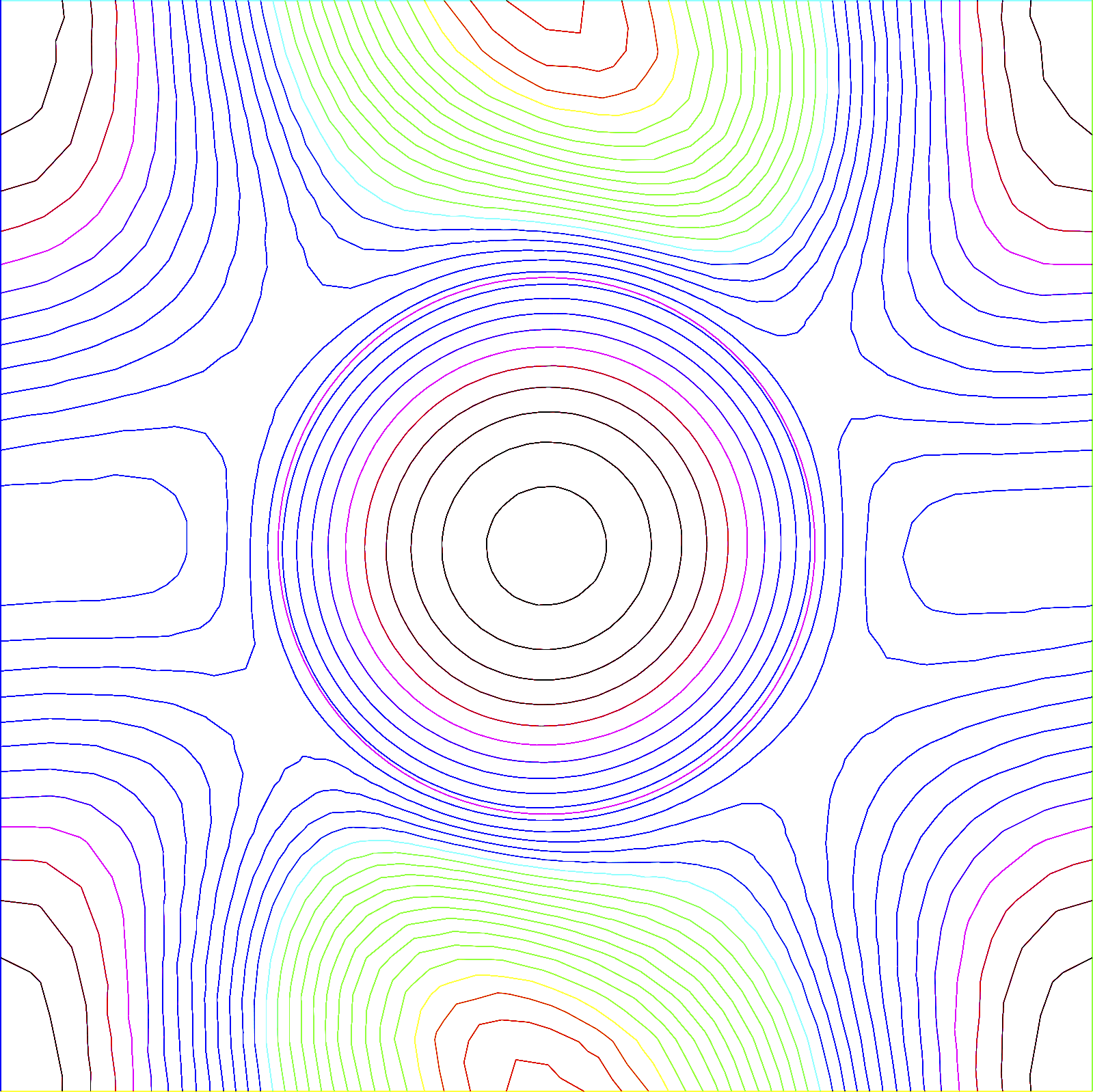}}
     {$\xi^{12}$}
&
\subf{\includegraphics[width=45mm]{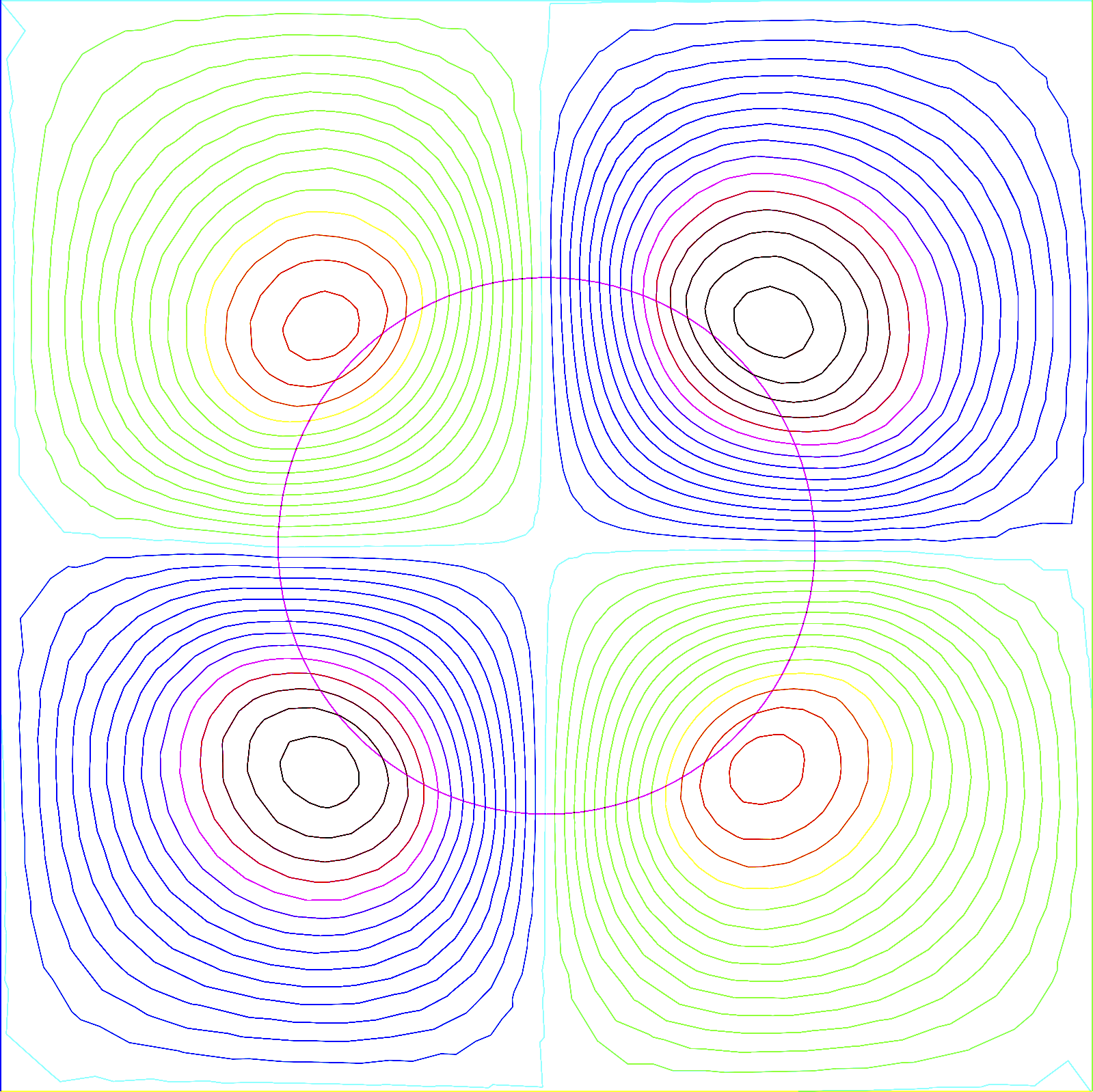}}
     {$\xi^{22}$}
\end{tabular}
\caption{Streamlines of the solution of the local problems for circular iron particles of $19\%$ volume fraction and $\alpha=1$ for a {\bf uniform distribution} of particles. The top row showcases the streamlines for the local solution $\vc{\chi^{m\ell}}$ in equation \eqref{var_chi} while the bottom row showcases the streamlines for the local solution $\vc{\xi^{m\ell}}$ in equation \eqref{var_xi}.}
\end{figure}

\begin{figure}[!htb]
\label{fig:tbh}
\centering
\begin{tabular}{ccc}
\subf{\includegraphics[width=63mm]{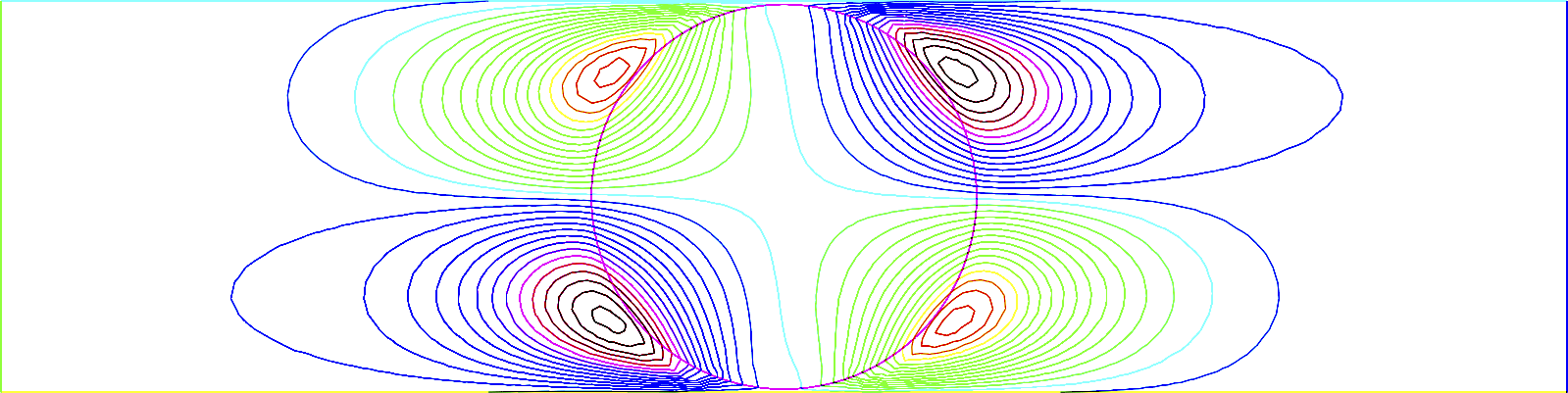}}
     {$\chi^{11}$}
&
\subf{\includegraphics[width=63mm]{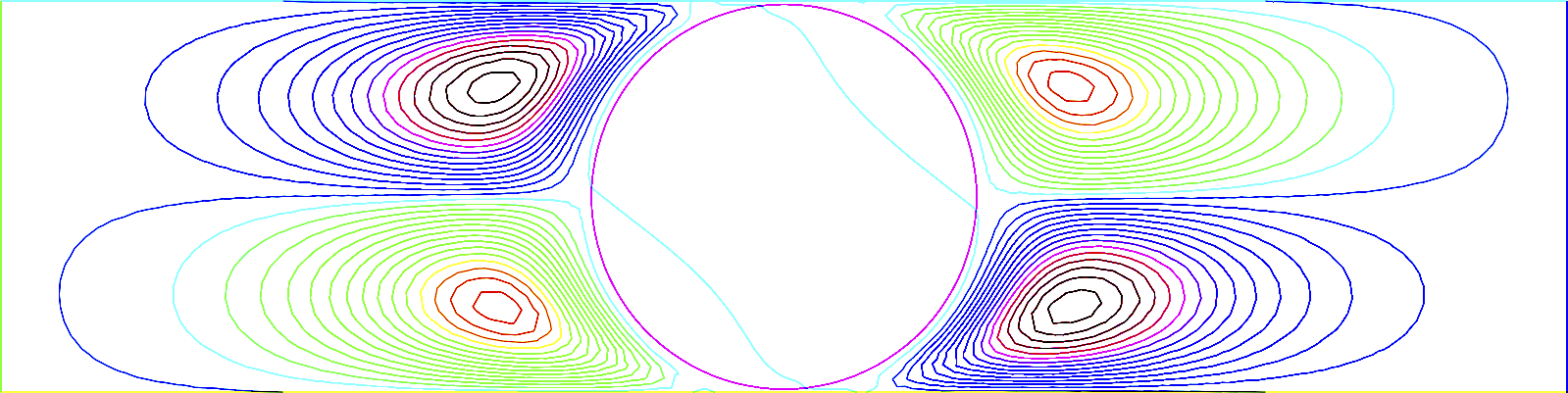}}
     {$\xi^{11}$}
\\
\subf{\includegraphics[width=63mm]{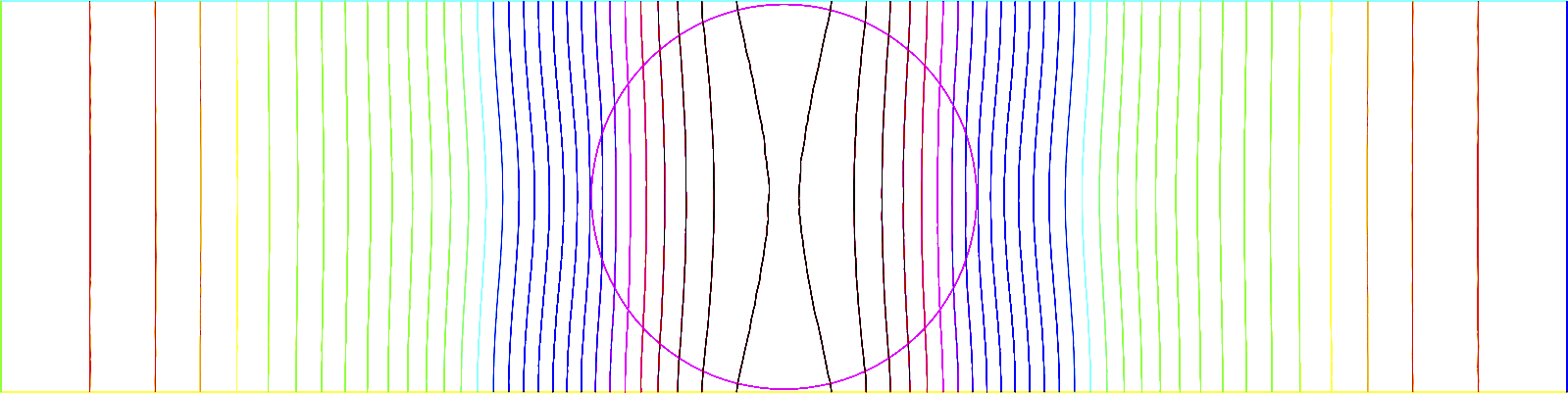}}
     {$\chi^{12}$}
&
\subf{\includegraphics[width=63mm]{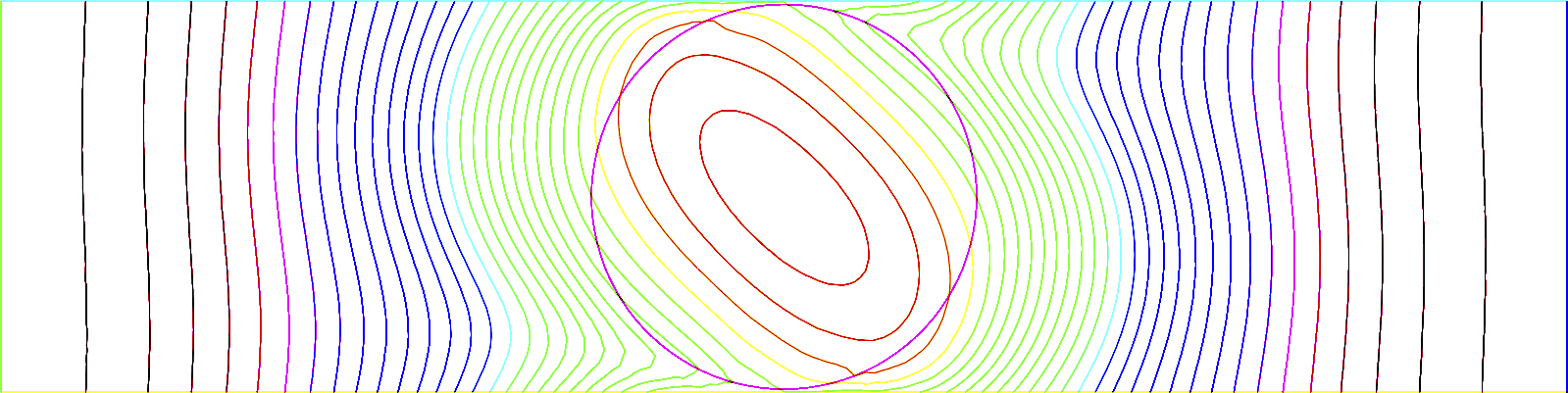}}
     {$\xi^{12}$}
\\
\subf{\includegraphics[width=63mm]{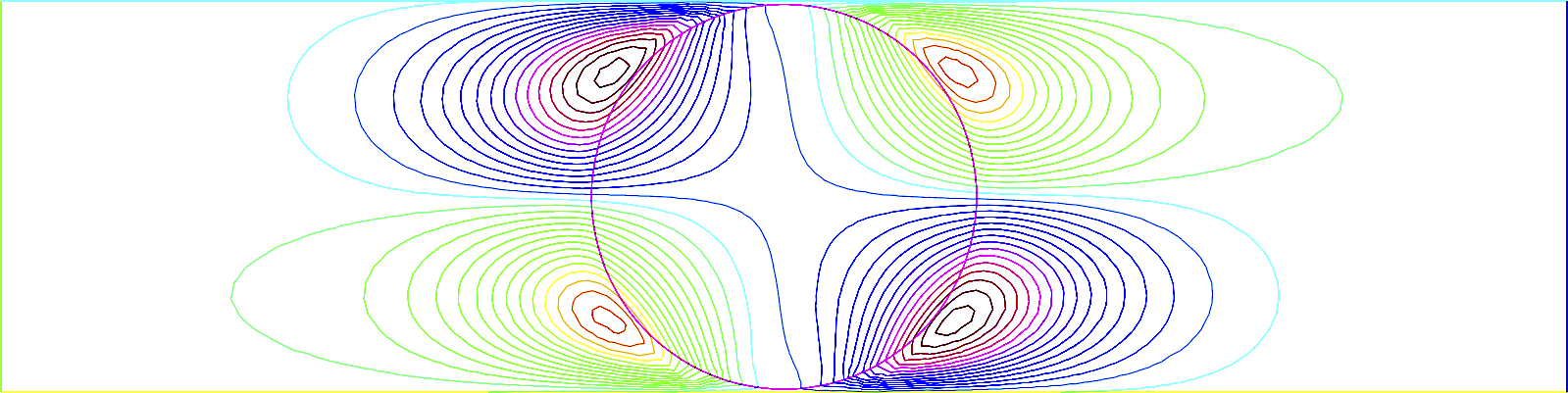}}
     {$\chi^{22}$}
&
\subf{\includegraphics[width=63mm]{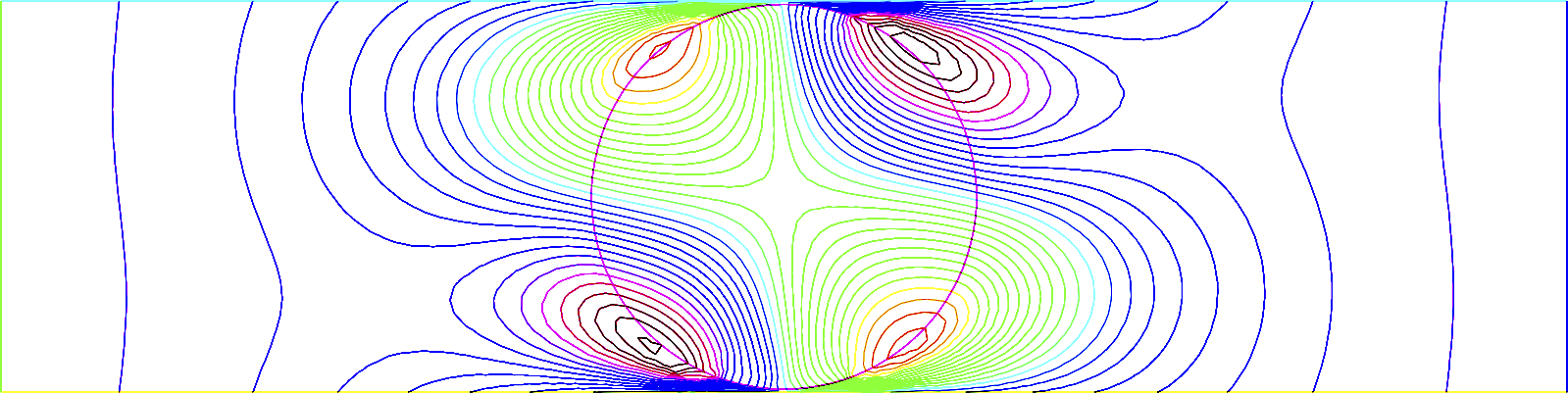}}
     {$\xi^{22}$}
\end{tabular}
\caption{Streamlines of the solution of the local problems for circular iron particles of $19\%$ volume fraction and $\alpha=1$ for particles distributed in {\bf chains}. The left row showcases the streamlines for the local solution $\vc{\chi^{m\ell}}$ in \eqref{var_chi} while the row on the right showcases the streamlines for the local solution $\vc{\xi^{m\ell}}$ in \eqref{var_xi}.}
\end{figure}

We compute the solutions of the local problems both for a uniform distribution of particles, Figure~\ref{fig:obo}, and for particles distributed in {\it chains}, Figure~\ref{fig:tbh}.  

\subsection{Influence of chain structures on the effective coefficients}

\begin{figure}[!htb]
\centering
\begin{tabular}{cc}
\subf{\includegraphics[width=69mm]{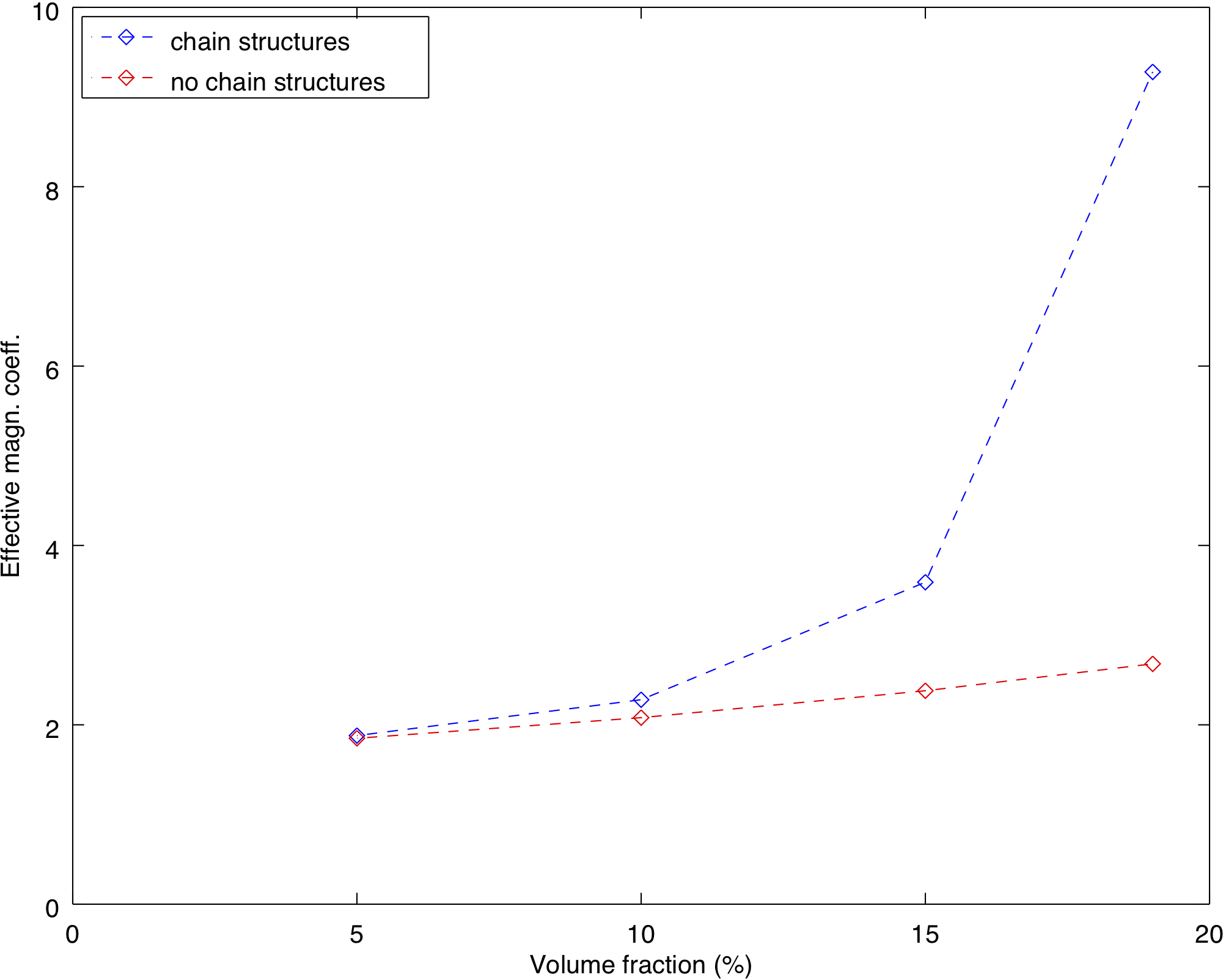}}
     {(a)}
&
\subf{\includegraphics[width=69mm]{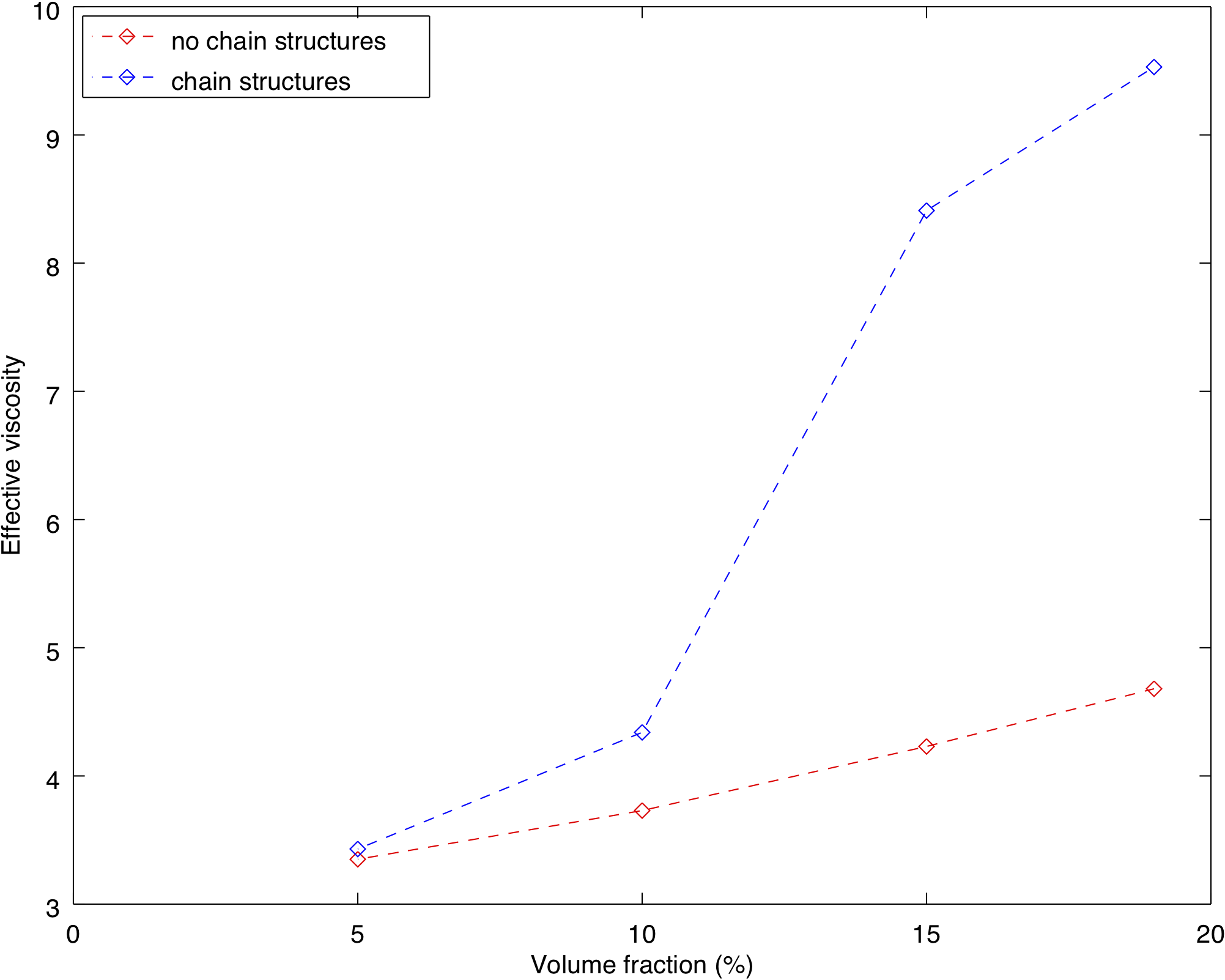}}
     {(b)}
\end{tabular}
\caption{Effective magnetic coefficient $\beta_s$ and effective viscosity $\nu_s$ plotted against volume fractions of $5\%$, $10\%$, $15\%$, and $19\%$ for circular iron particles. On image (a), the {\bf red} color curve showcases the increase of the effective magnetic coefficient $\beta_s$ under uniform particle distribution while the {\bf blue} color curve showcases the increase of the effective magnetic coefficient $\beta_s$ in the presence of {\bf chain structures}. Similarly, on image (b) the {\bf red} color curve showcases the increase of the effective viscosity $\nu_s$ under uniform particle distribution while the {\bf blue} color curve showcases the increase of the effective viscosity $\nu_s$ in the presence of {\bf chain structures}.}
\label{fig:coeff}
\end{figure}

Here we are interested in exploring the influence {\it chain structures} have on the effective coefficients. In Figure~\ref{fig:coeff} we plot the effective coefficients $\beta_s$ and $\nu_s$ versus different volume fractions for the two types of particle distributions. 

In Figure~\ref{fig:coeff} we see that the presence of {\it chain structures} rapidly increase the effective viscosity $\nu_s$ and the effective magnetic coefficient $\beta_s$ as the volume fraction increases in a non-linear fashion. This is in stark contrast to the case where no particle {\it chains} are present where we notice a linear increase of the effective coefficients. For a particle of $19\%$ volume fraction the presence of {\it chain structures} leads to roughly a $104\%$ increase of the effective viscosity $\nu_s$ and a $246\%$ increase of the magnetic coefficient $\beta_s$.  

\subsubsection{Velocity profiles of Poiseuille and Couette flows}

In this section we compute the cross sectional velocity profiles of Poiseuille and Couette flow for circular suspensions of rigid particles. We denote by $\vc{v} = (v_1,v_2)$ the two dimensional velocity and by $\vc{H} = (H_1,H_2)$ the two dimensional magnetic field. Hence, the two dimensional stress of equation~\eqref{eq:vecform} reduces to,
\[
	\sigma^H + \tau^H = -\bar{p}^0 \, I + \nu_s \, e(\vc{v^0}) - \frac{1}{2} \, \beta_s \,  \left|\widetilde{\vc{H}}^0\right|^2 \, I + \beta_s \, \widetilde{\vc{H}}^0 \otimes \widetilde{\vc{H}}^0
\]
and the two dimensional magnetorheological equations in \eqref{eq:vecform} reduce to the following:

\begin{subequations} 
\label{eq:MR2D}
\begin{align} 
\frac{\nu_s}{2} \, \left( \frac{\partial^2 v_1}{\partial x_1^2} + \frac{\partial^2 v_1}{\partial x_2^2} \right) -\frac{\partial \, \pi^0}{\partial x_1}  + \frac{\beta_s}{2} \, \frac{\partial ( H_1^2 - H_2^2) }{\partial x_1} + \beta_s \, \frac{\partial (H_1 \, H_2) }{\partial x_2} &= 0, \label{FE1}\\
\frac{\nu_s}{2} \, \left( \frac{\partial^2 v_2}{\partial x_1^2} + \frac{\partial^2 v_2}{\partial x_2^2} \right) -\frac{\partial \, \pi^0}{\partial x_2} + \beta_s \, \frac{\partial (H_1 \, H_2)}{\partial x_1} + \frac{\beta_s}{2} \, \frac{\partial ( H_2^2 - H_1^2 )}{\partial x_2} &= 0, \label{FE2}\\
\frac{\partial}{\partial x_1} \left ( \mu^H \, H_1 \right ) + \frac{\partial}{\partial x_2} \left ( \mu^H \, H_2 \right ) &= 0, \label{MD}\\
\frac{\partial H_2}{\partial x_1} - \frac{\partial H_1}{\partial x_2} = R_m \, (\mu^{HS}_{22} \, v_1 \, H_2 - \mu^{HS}_{11} \, v_2& \, H_1), \label{MC}\\
\frac{\partial v_1}{\partial x_1} + \frac{\partial v_2}{\partial x_2} &= 0. \label{FD}
\end{align}
\end{subequations}      

\subsubsection{Poiseuille flow}

We consider the problem of a steady flow due to a pressure gradient between two infinite, parallel, stationary plates that are non-conducting and non-magnetizable with one plate aligned along the $x_1$--axis while the second plate is of distance one unit apart. We apply a stationary magnetic field $\vc{H}$ on the bottom plate. Since we are dealing with infinite plates, the velocity $\vc{v}$ depends only on $x_2$. Using \eqref{FD} we immediately obtain that $v_2$ is constant and since the plates are stationary $v_2=0$. Since the flow is unidirectional, we expect that the the magnetic field will depend only on the height $x_2.$ Hence, using \eqref{MD} we obtain $H_2 (x_2) = K$, while the component parallel to the flow depends on the fluid velocity. Therefore the equations in \eqref{eq:MR2D} reduce to the following, 
\begin{subequations} \label{eq:MR2DPC}
\begin{align} 
 \frac{\nu_s}{2} \, \frac{\partial^2 v_1}{\partial x_1^2} + \beta_s \, K \, \frac{\partial H_1}{\partial x_2} &= \frac{\partial \, \pi^0}{\partial x_1}, \label{FE1_PC} \\
-\frac{\partial \, \pi^0}{\partial x_2} - \frac{1}{2} \, \beta_s \, \frac{\partial H_1^2}{\partial x_2} &= 0, \label{FE2_PC} \\
-\frac{\partial H_1}{\partial x_2} &= R_m \, \mu^{HS}_{22} \, K \, v_1. \label{MC_PC}
\end{align}
\end{subequations} 

Making use of \eqref{FE2_PC} we obtain that $\pi^0 (x_1,x_2) + \frac{1}{2} \, \beta_s \, H_1(x_2)^2$ is a function of only $x_1$ and therefore by differentiating the expression with respect to $x_1$ we get that $\frac{\partial \, \pi^0}{\partial x_1}$ is a function only $x_1$. Therefore, on \eqref{FE1_PC}
the left hand side is a function of $x_2$ and the right hand side is a function of $x_1$. Thus, they have to be constant. Substituting \eqref{MC_PC} in \eqref{FE1_PC} we obtain the following differential equations,
\begin{equation} \label{s_ode} 
 \frac{d^2 \, v_1}{d \, x_2^2} - \lambda^2 \, v_1 = C_p, \quad \frac{\partial \, \pi^0}{\partial x_1} = C_p \text{ with } \lambda = \sqrt{\frac{2 \, R_m \, \mu^{HS}_{22} \, \beta_s}{\nu_s}} \, K
\end{equation}

The general solution of \eqref{s_ode} is, 
\[ 
v_1(x_2) = c_1 \, e^{\lambda \, x_2} + c_2 \, e^{-\lambda \, x_2} + \frac{C_p}{\nu_s \, \lambda^2}. 
\] 

Given that $v_1(0)=v_1(1) = 0$ we have,
\begin{equation} \label{Pois}
v_1(x_2) = \frac{C_p}{\nu_s \, \lambda^2} \left(\frac{\sinh(\lambda \, x_2) - \sinh(\lambda \, (x_2 - 1)) }{\sinh(\lambda)} - 1 \right). 
\end{equation} 

\begin{remark}
As $\lambda$ tends to zero we have $\lim_{\lambda \to 0} v_1(x_2) = \frac{ C_p }{2 \, \nu_s} \, x_2 \, (x_2 - 1)$, which is precisely the profile of Poiseuille flow. 
\end{remark}

Once the velocity $v_1(x_2) $ is known, we can use \eqref{MC_PC} to compute $H_1(x_2)$ with boundary condition $H_1(0) = K_1$ and obtain,
\[ 
H_1(x_2) = R_m \, \mu^{HS}_{22} \, K \, \frac{C_p}{\nu_s \, \lambda^3 \, \sinh(\lambda)} \, (-\cosh(\lambda \, x_2) + \cosh(\lambda \, (x_2 - 1)) - \cosh(\lambda) + 1) +  K_1.
\]

\subsubsection{Couette flow}
The setting and calculations for the unidirectional Couette flow are the same as Poiseuille flow. In a similar way, we can carry out computations for the plane Couette flow. For simplicity we assumed the bottom plate is the $x_1$ axis and the top plate is at $x_2=1$ and the pressure gradient is zero. A shear stress $\gamma$ is applied to the top plate while the bottom plate remains fixed. Thus, we solve \eqref{s_ode} with initial 
conditions $v_1(0)=0$ and $v'_1(1) = \gamma$ and obtain
\begin{equation}\label{Couette}
v_1(x_2) = - \frac{(C_p e^{-\lambda} - \gamma \nu_s)e^{\lambda x_2}}{\lambda(e^{\lambda} + e^{-\lambda})\nu_s} - \frac{(C_p e^{\lambda} + \gamma \nu_s)e^{-\lambda x_2}}{\lambda(e^{\lambda} + e^{-\lambda})\nu_s} + \frac{C_p}{\nu_s \lambda^2}
%\frac{\dot{\gamma} \, \nu_s \, \lambda \, \sinh(\lambda x_2) + C_p \, \cosh(\lambda \, (x_2 - 1) )}{\nu_s \lambda^2 \cosh(\lambda)} - \frac{C_p}{\nu_s \, \lambda^2} 
\end{equation}

\begin{remark}
We remark that for a zero pressure gradient the limit as $\lambda$ approaches zero the limit of $v_1(x_2) = \gamma \, x_2$, the profile of regular Couette flow.
\end{remark}

To compute $H_1$ we use \eqref{MC_PC} to obtain,
\begin{gather}
\begin{aligned} 
H_1(x_2) (R_m \mu^{HS}_{22} K)^{-1} 
& = \frac{(C_p e^{-\lambda} - \gamma \nu_s)e^{\lambda x_2}}{\lambda^2 (e^{\lambda} + e^{-\lambda})\nu_s} - \frac{(C_p e^{\lambda} + \gamma \nu_s)e^{-\lambda x_2}}{\lambda^2 (e^{\lambda} + e^{-\lambda})\nu_s} - \frac{C_p x_2}{\nu_s \lambda^2} \\
& - \frac{(C_p e^{-\lambda} - \gamma \nu_s)}{\lambda^2 (e^{\lambda} + e^{-\lambda})\nu_s} + \frac{(C_p e^{\lambda} + \gamma \nu_s)}{\lambda^2 (e^{\lambda} + e^{-\lambda})\nu_s} + K_1 (R_m \mu^{HS}_{22} K)^{-1} 
%\frac{\dot{\gamma}}{\lambda^2 \, \cosh(\lambda)} \, (\cosh(\lambda \, x_2) - 1) 
%+ \frac{C_p}{\nu_s \, \lambda^3 \, \cosh(\lambda)} \, (\sinh(\lambda \, (x_2 - 1)) - \sinh(\lambda)) \nonumber \\
%& - \frac{C_p \, x_2}{\nu_s \, \lambda^2} + K_1. \nonumber
\end{aligned}
\end{gather}

We plot the velocity profiles for Poiseuille and Couette flows, computed in the previous paragraph, for a magnetorheological suspension of iron particles for different magnetic field intensities. The volume fraction of the particles is set to $19\%$. Carrying out explicit computations of the effective coefficients in \eqref{nu}, \eqref{beta} and \eqref{mu:s} we obtain,  
\vspace{.1cm}
\begin{center}
\begin{tabular}{ l | c  c  c  c }
                             ~    & $\nu_s $ & $\beta_s $ & $\mu^{HS}_{22}$ \\ \hline
  $2 \times 1/2$ unit cell & 9.53           & 9.28                   & 3.71 \\
  $1 \times 1$ unit cell    & 4.68           & 2.66                   & 0.47 \\
\end{tabular}
\end{center}
\vspace{.3cm}

\begin{figure}[!htb]
\centering
\begin{tabular}{ccc}
\subf{\includegraphics[width=45mm, height=45mm]{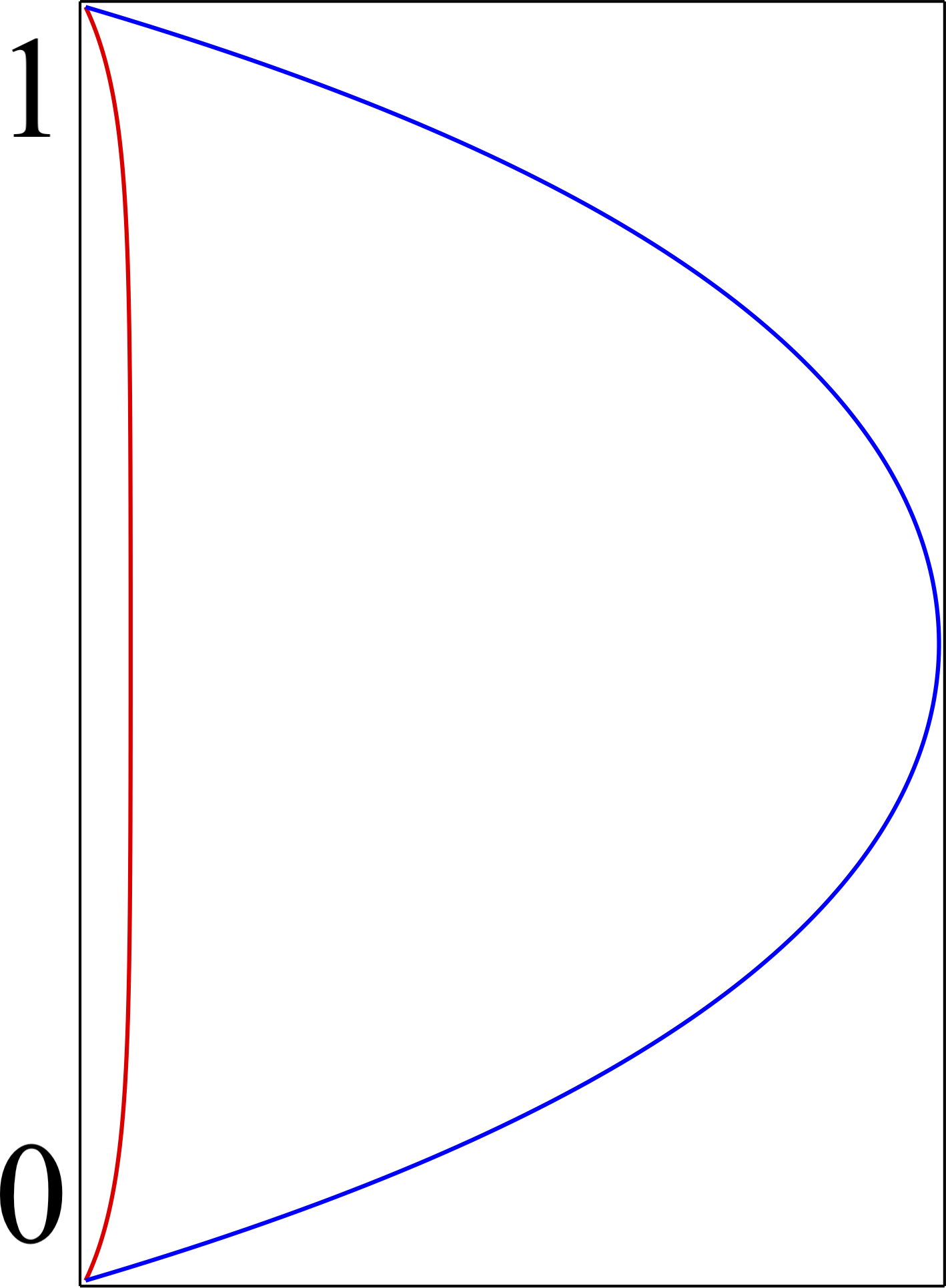}}
     {$H_2 = 50$}
&
\subf{\includegraphics[width=45mm, height=45mm]{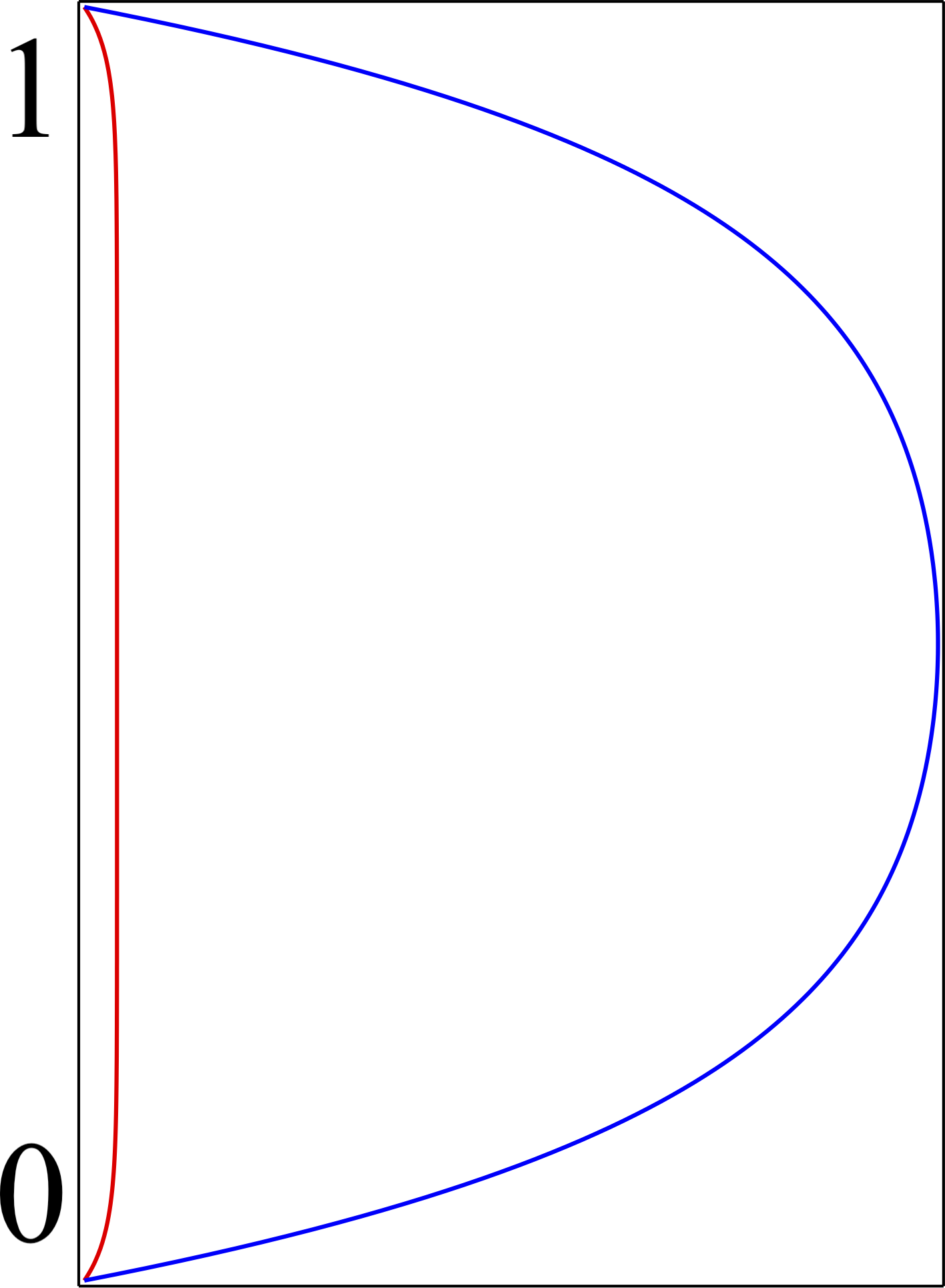}}
     {$H_2 = 100$}
&
\subf{\includegraphics[width=45mm, height=45mm]{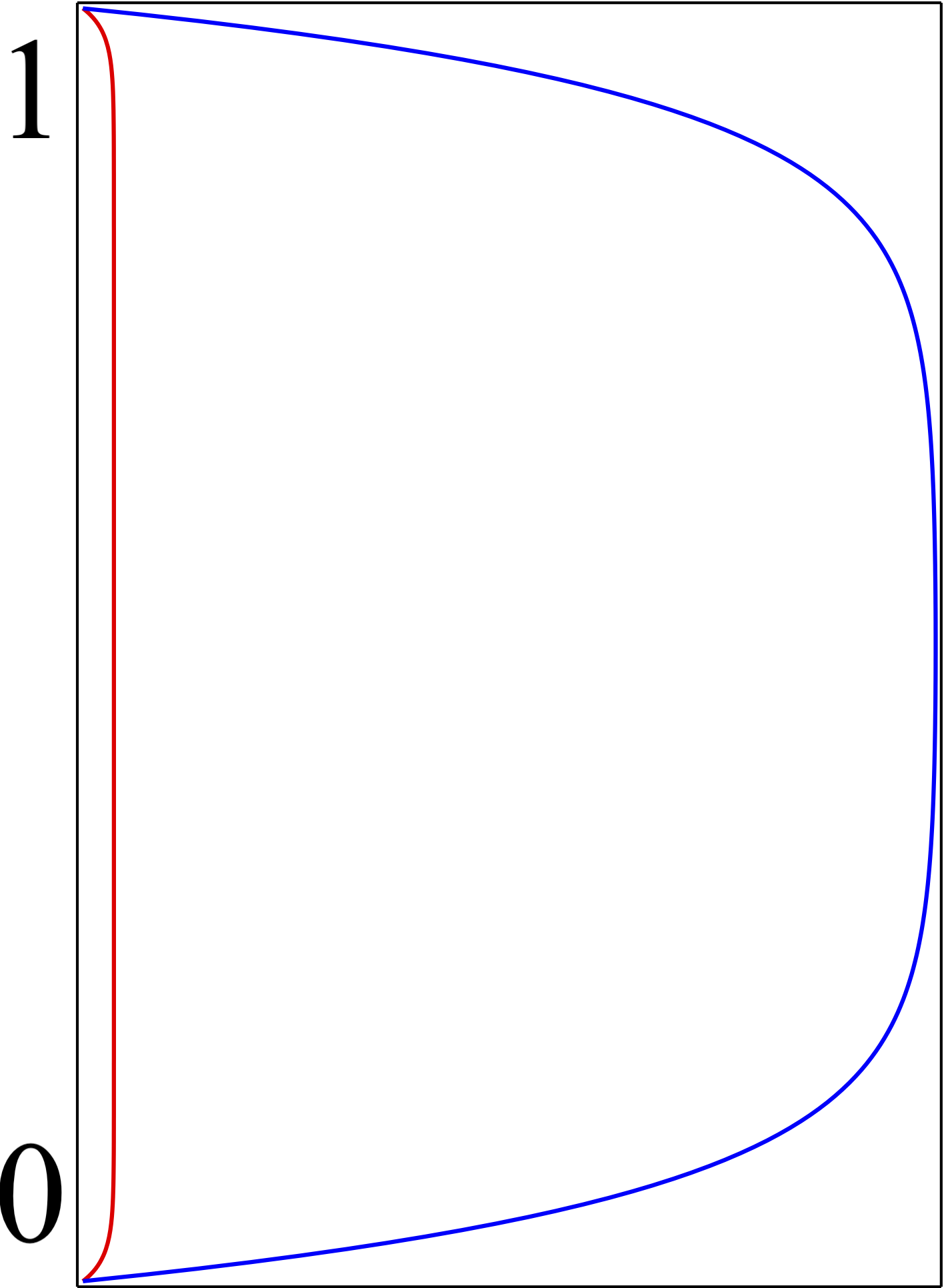}}
     {$H_2 = 200$}
\end{tabular}
\caption{Velocity profile of {\bf Poiseuille flow} for a magnetorheological suspension of circular iron particles of $19\%$ volume fraction in a viscous non-conducting fluid with $\alpha=1$ and $R_m=10^{-2}$ plotted against different values of the magnetic field. The {\bf blue} color curve represents the velocity profile for a uniform particle distribution while the {\bf red} color curve represents the velocity profile for a particle distribution in {\bf chains}}.
\label{fig:PF}
\end{figure}

Figure~\ref{fig:PF} showcases the velocity profile of magnetorheological Poiseuille flow for a constant pressure gradient and for three different values of the magnetic field. We can see that the damping force increases as we increase the magnetic field; the profile is close to flat in the middle region for high values of the magnetic field, but it is not parabolic close to the walls as in the case of Bingham flows. Moreover, it seems that the presence of {\it chain structures} (red color) turns the magnetorheological fluid into a stiffer gel-like structure~\cite{Bossis} at lower intensity magnetic fields. This phenomenon results in a Poiseuille flow that is much slower in the presence of {\it chain structures}.

\begin{figure}[!htb]
\centering
\begin{tabular}{ccc}
\subf{\includegraphics[width=45mm]{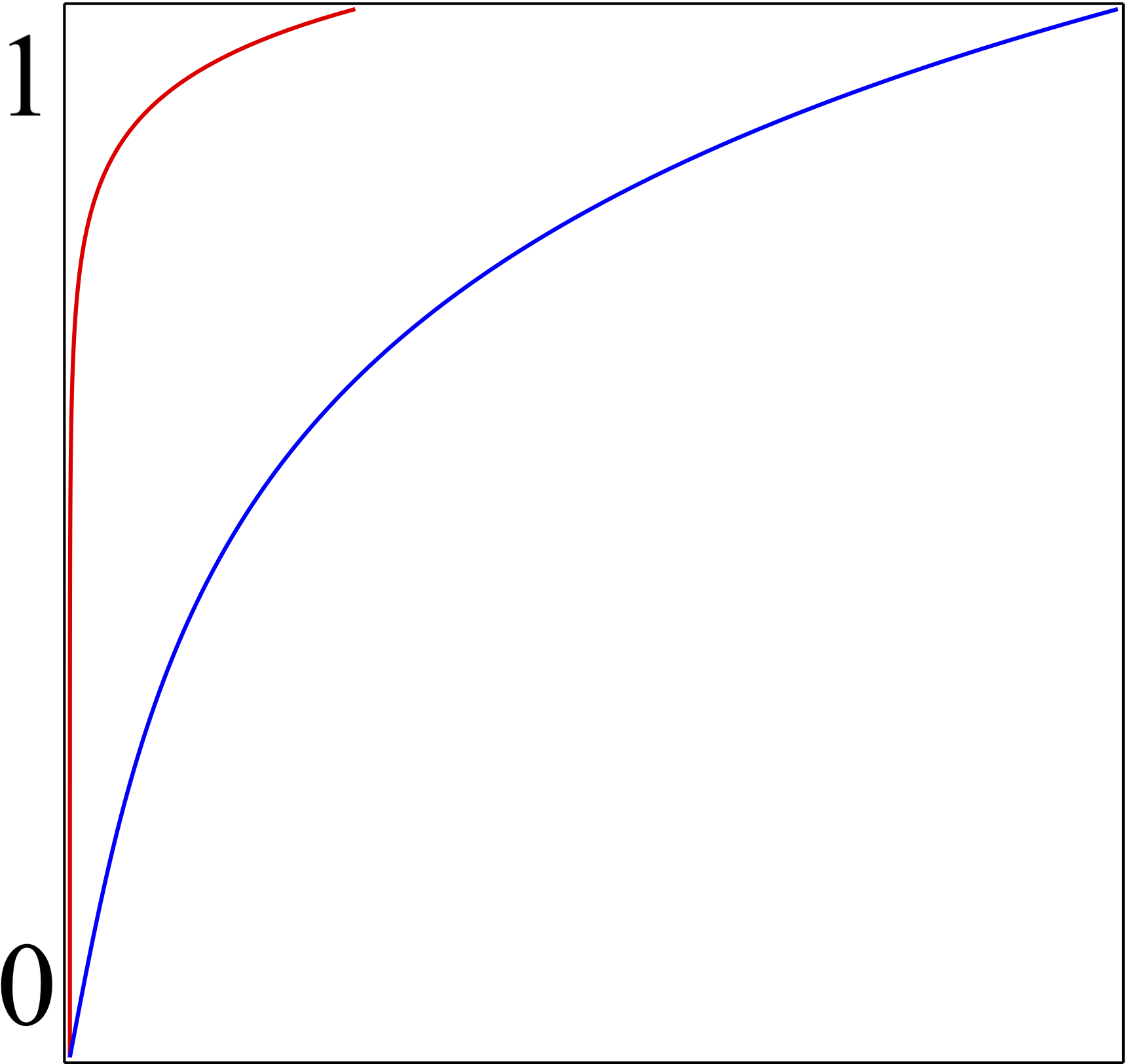}}
     {$H_2 = 50$}
&
\subf{\includegraphics[width=45mm]{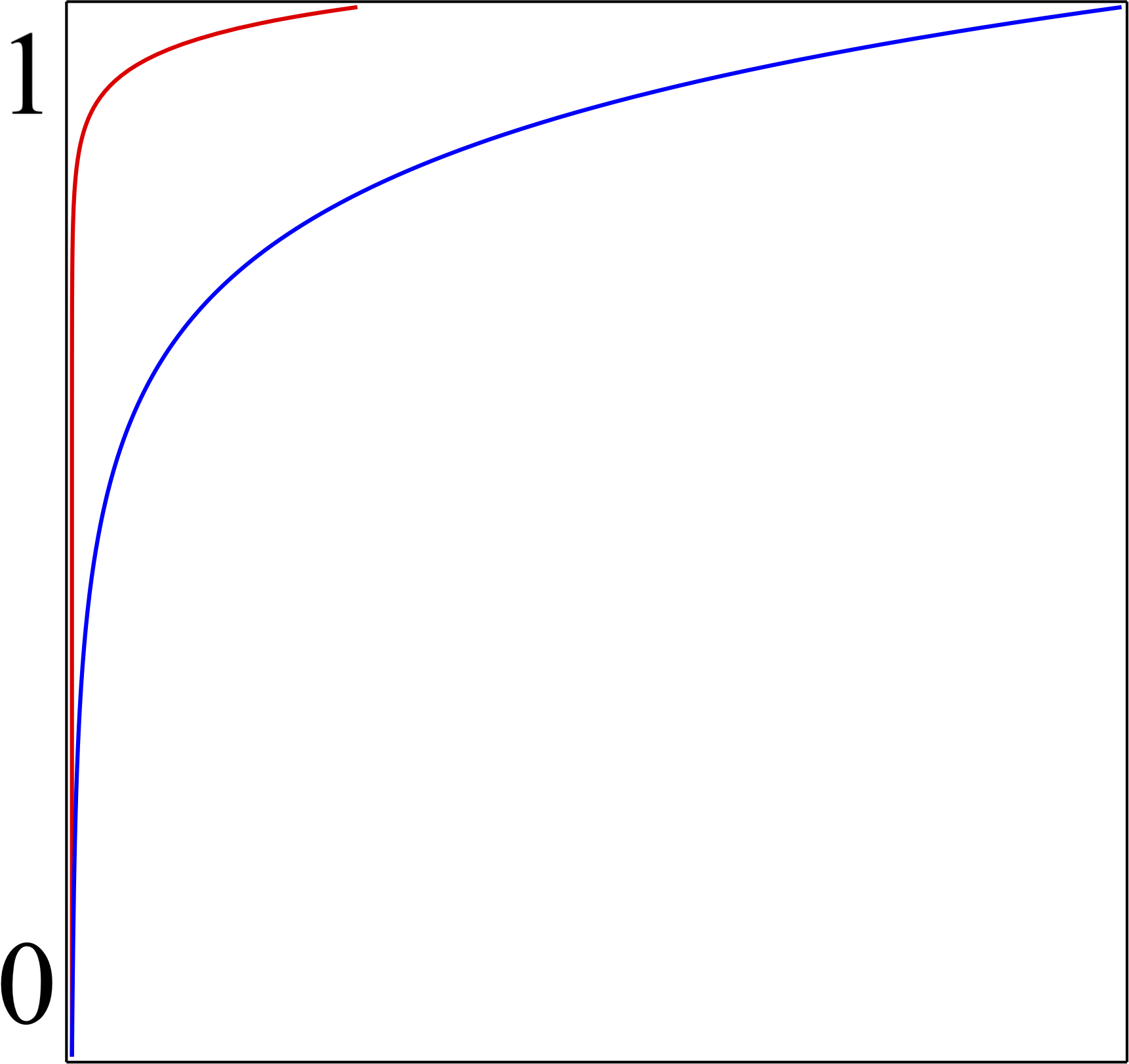}}
     {$H_2 = 100$}
&
\subf{\includegraphics[width=45mm]{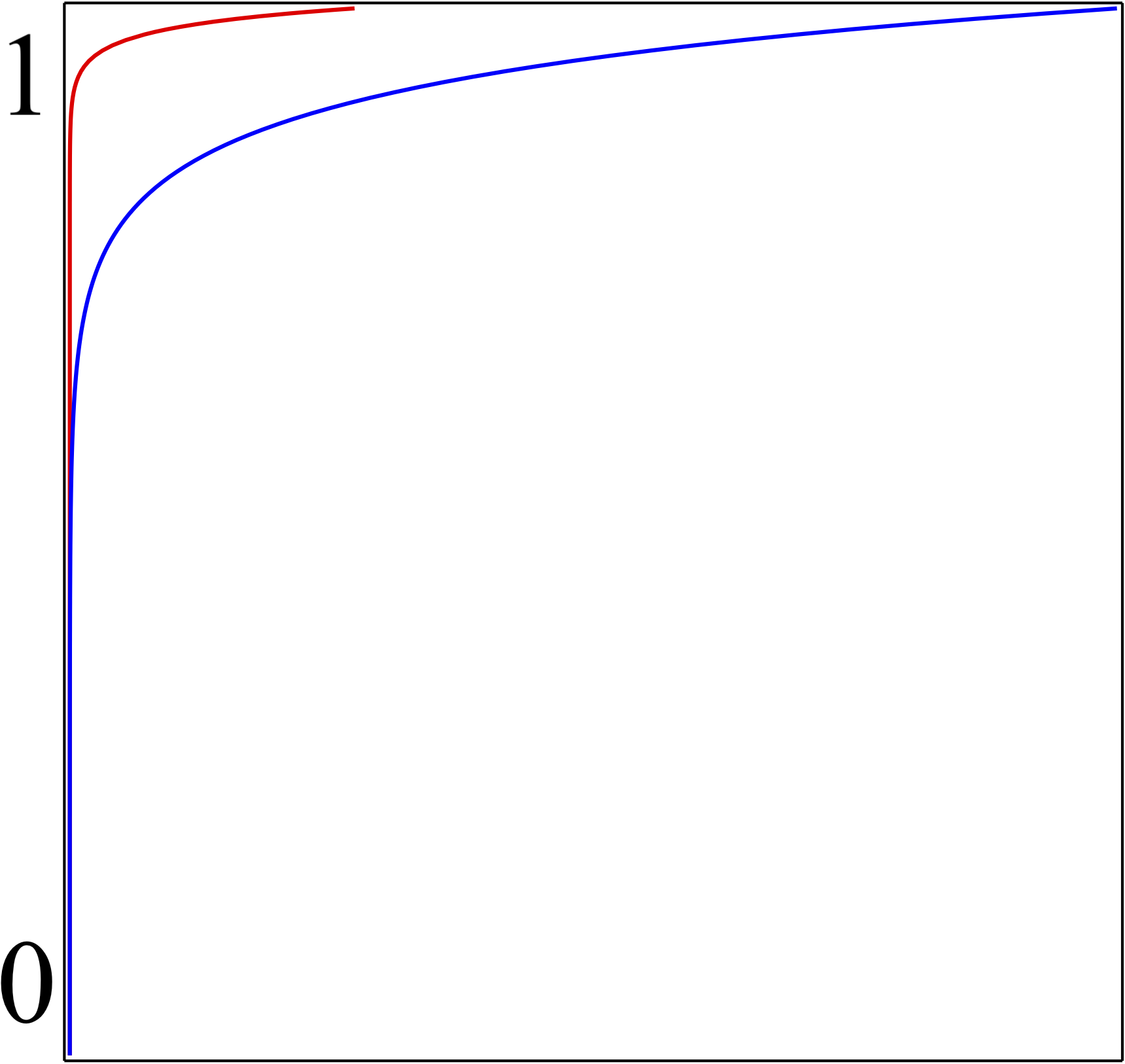}}
     {$H_2 = 200$}
\end{tabular}
\caption{Velocity profile of {\bf Couette flow} for a magnetorheological suspension of iron particles of $19\%$ volume fraction in a viscous non-conducting fluid with $\alpha=1$ and $R_m=10^{-2}$ plotted against different values of the magnetic field. The {\bf blue} color curve represents the velocity profile for a uniform particle distribution while the {\bf red} color curve represents the velocity profile for a particle distribution of {\bf chains}}. 
\label{fig:CF}
\end{figure}

Figure~\ref{fig:CF} showcases the velocity profile of magnetorheological Couette flow for zero pressure gradient and three different values of the magnetic field. We can observe that as the magnetic filed increases the magnetorheological fluid becomes harder to shear. Hence, an {\it apparent yield stress} is present. The {\it apparent yield stress} is larger in the presence of {\it chain structures} (red color) than in the absence of {\it chain structures} (blue color). However, the velocity profile is not linear like in the case of Bingham fluids.

\begin{figure}[!htb]
\centering
\begin{tabular}{cc}
\subf{\includegraphics[width=62mm]{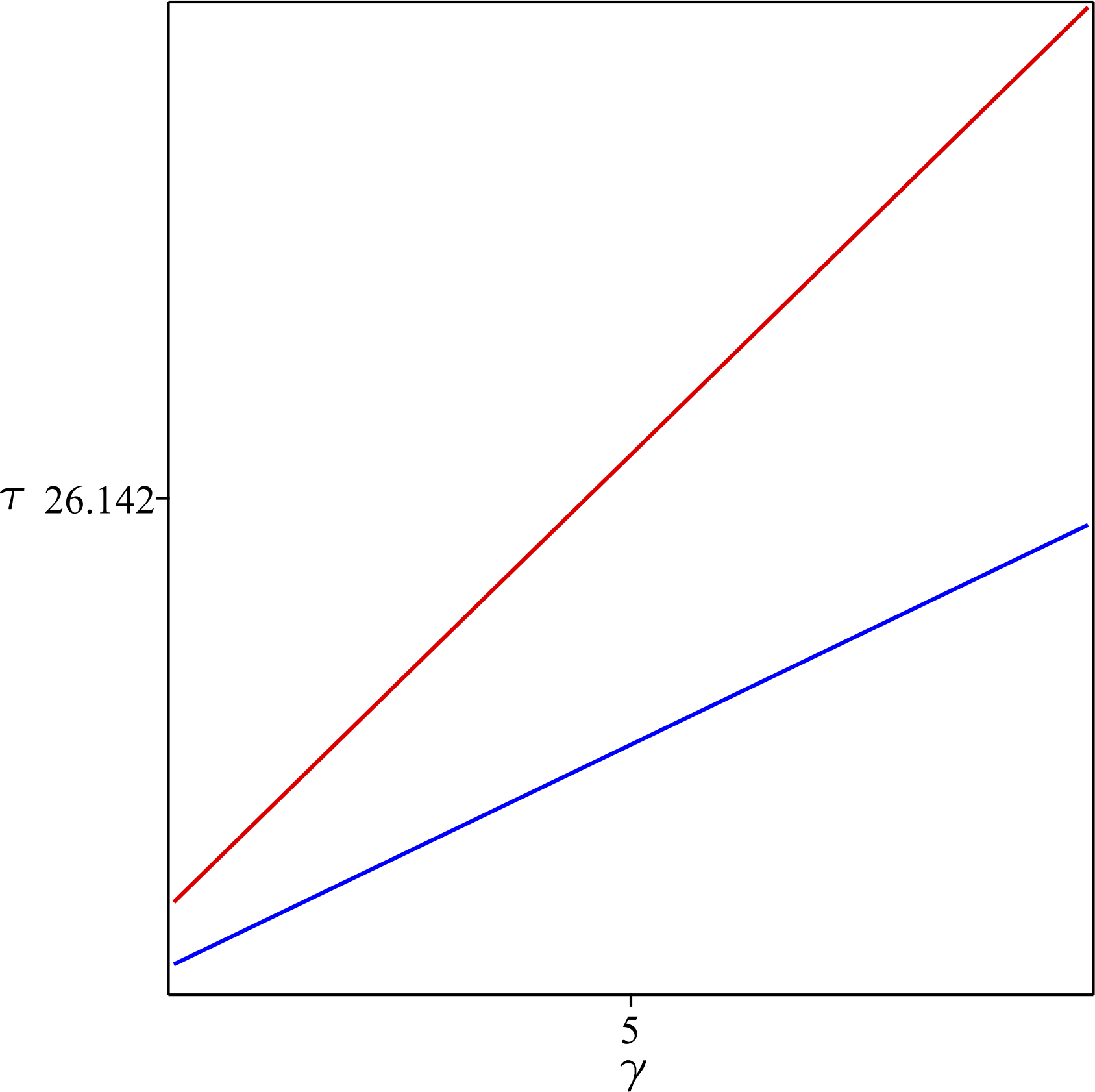}}
     {$H_2 = 50$}
&
\subf{\includegraphics[width=62mm]{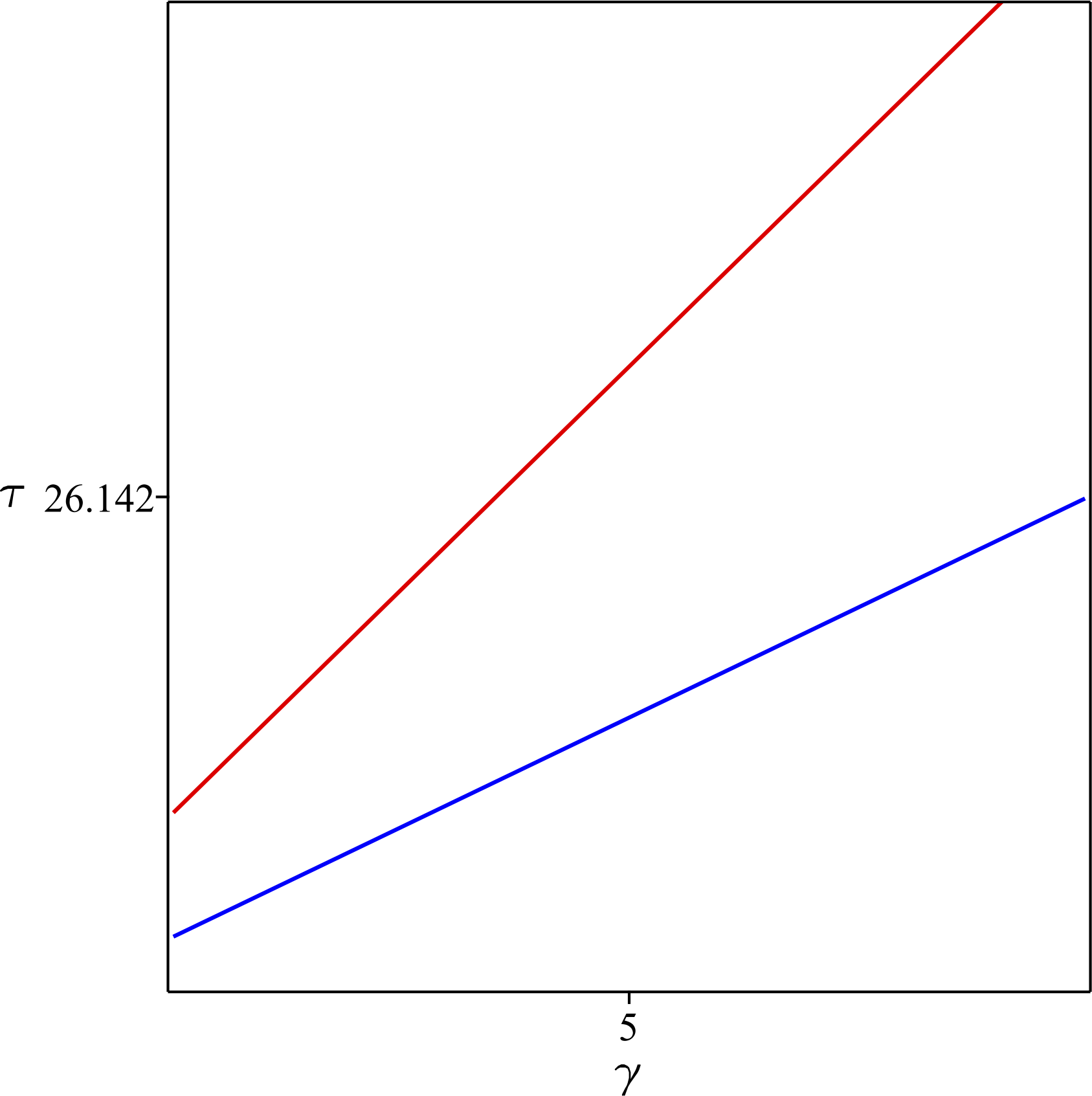}}
     {$H_2 = 100 $} \\
\subf{\includegraphics[width=62mm]{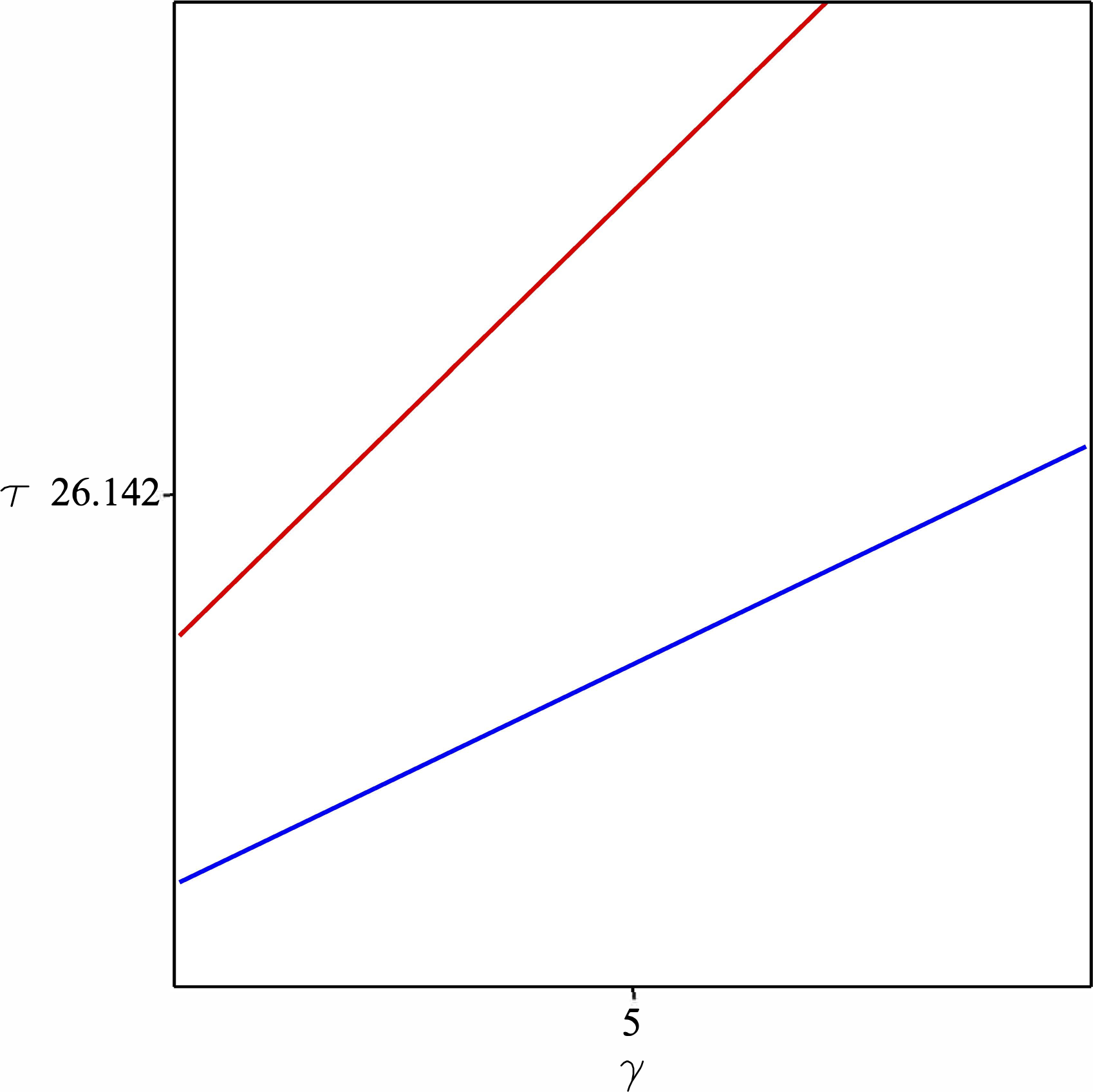}}
     {$H_2 = 200$}
&
\subf{\includegraphics[width=63mm, height=62mm]{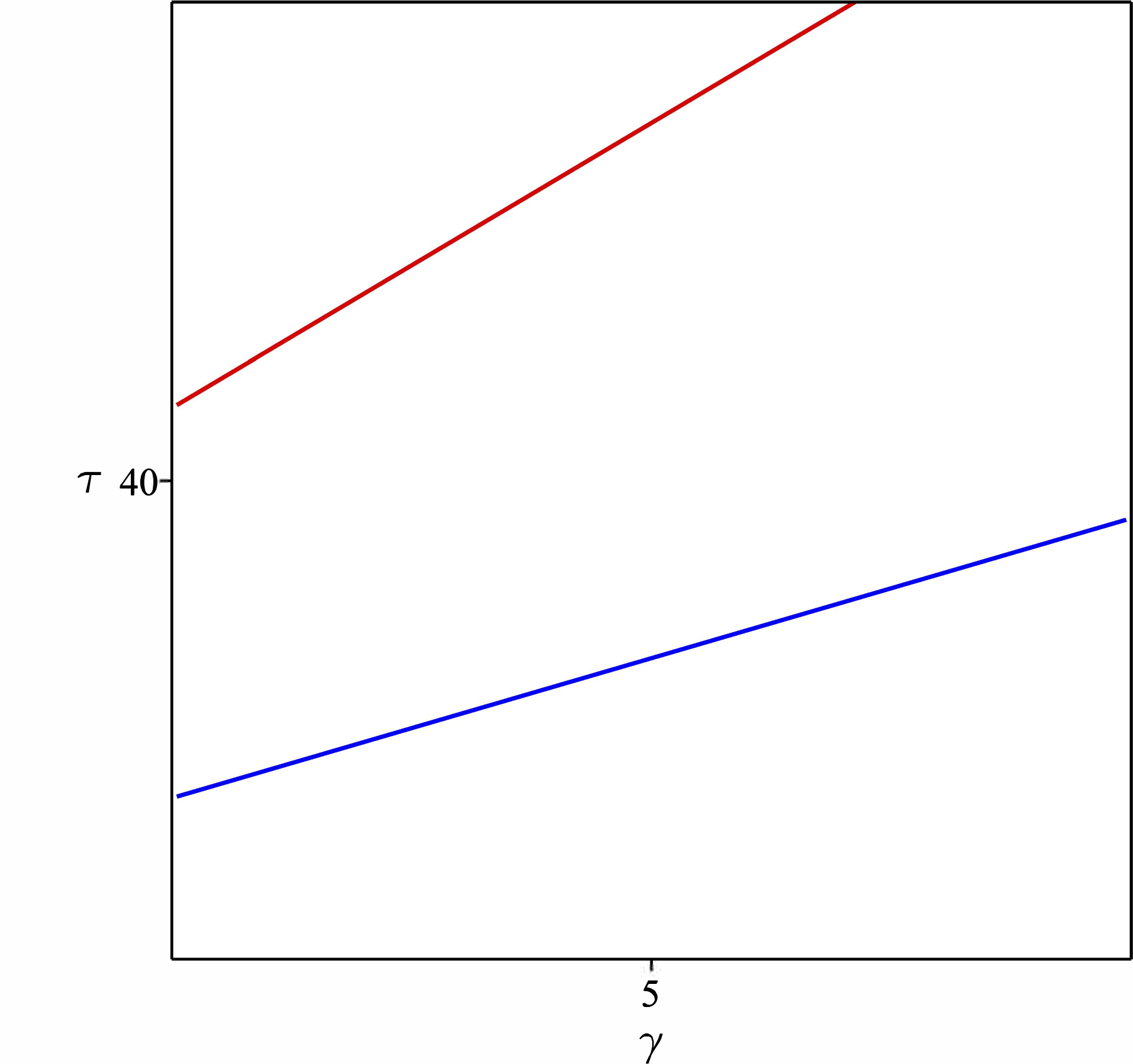}}
     {$H_2 = 500$}
\end{tabular}
\caption{{\bf Shear stress} $\tau$ versus {\bf shear rate} $\gamma$ for a magnetorheological suspension of circular iron particles of $19\%$ volume fraction in a viscous non-conducting fluid with $\alpha=1$, $R_m=10^{-2}$, and $K_1=10^{-2}$. The {\bf red} color curve represents particle distribution in {\bf chains} while the {\bf blue} color curve represents a uniform particle distribution. One can observe that the {\bf yield stress} of the magnetorheological fluid is consistently higher in the presence of {\bf chain structures} (red color).}
\label{fig:SR}
\end{figure}

Figure.~\ref{fig:SR} showcases the {\it shear stress} $\tau$ vs {\it shear rate} $\gamma$ relationship measured at $x_2=1$. When $K_1=0$ there is no yield stress present. However, for very small non-zero values of $K_1$ we obtain results similar to~\cite{YLL05} for the linear portion of the shear stress versus shear rate curve. Additionally, we can see that the presence of {\it chain structures} (red color) produce a higher yield stress as well as a steeper slope that gives higher shear stress values as the shear rate increases. 

For shear experiments, the response of magneto-rheological fluids is often modeled using a Bingham constitutive law~\cite{Bossis}, \cite{GS15}, \cite{Ha92}, \cite{PV02}. Although the Bingham constitutive law measures the response of the magnetorheological fluid quite reasonably, actual magnetorheological fluid behavior exhibits departures from the Bingham model \cite{HL11}, \cite{YLL05}. In Figure~\ref{fig:PF} and Figure~\ref{fig:CF} we see that for low intensity magnetic fields the Bingham constitutive law is not adequate, however, it appears that for higher intensity magnetic fields the flow gets closer to resembling a Bingham fluid behavior.

\section{Conclusions}

We considered a suspension of magnetizable iron particles in a non-magnetizable, non-conducting aqueous viscous fluid. We obtained the effective equations governing the behavior of the magnetorheological fluid presented in equation~\eqref{eq:vecform}. 
The material parameters can be computed from the local problems~\eqref{Maxwell_local}, \eqref{var_chi}, \eqref{var_xi}, derived from the balance of mass, momentum, and Maxwell equations. 

The proposed model generalizes the model put forth in~\cite{NR64}, \cite{Ro14} in two ways: 

\begin{itemize}
\item[$\bullet$] First by providing exact formulas for the effective coefficients which consist of the homogenized viscosity, $\nu^H$, and three homogenized magnetic permeabilities, $\mu^H$, $\mu^{HS}$, and $\beta^H$, which all depend on the geometry of the suspension, the volume fraction, the magnetic permeability $\mu$, the Alfven number $\alpha$, and the particles distribution. 
\item[$\bullet$] Second, by coupling the fluid velocity field with the magnetic field through Ohm's law. 
\end{itemize}

Using the finite element method we carried out explicit computations of the effective coefficients for spherical iron particles of different volume fractions both under uniform particle distribution and particle distribution in {\it chains} and showcased the nonlinear effect particle {\it chain structures} have in the effective coefficients as the volume fraction increases. 

In the case of {\it Poiseuille} flows we calculated the velocity profile explicitly; for small intensity magnetic fields the velocity profile is close to parabolic while for large intensity magnetic fields an {\it apparent yield stress} is present and the flow profile approaches a {\it Bingham flow} profile. The magnetorheological effect is significantly higher when {\it chain structures} are present. A similar analysis has been done for the {\it Couette} flows.

\section*{Acknowledgments}

G.N. was partially funded by the Deutsche Forschungsgemeinschaft (DFG, German Research Foundation) under Germany's Excellence Strategy – The Berlin Mathematics Research Center MATH+ (EXC-2046/1, project ID: 390685689) and would like to express his gratitude to Konstantinos Danas and Andrei Constantinescu for their fruitful discussions and suggestions.

\bigskip 

\bibliography{ref}

\end{document}